\documentclass [12pt] {report}
\usepackage{amssymb}
\newcommand {\D}[2] {\displaystyle\frac{\partial{#1}}{\partial{#2}}}

\newcommand {\Dd}[3] {\displaystyle\frac{\partial^2{#1}}{\partial{#2}\partial{#3}}}
\newcommand {\al} {\alpha}

\newcommand {\ga} {\gamma}
\newcommand {\la} {\lambda}

\newcommand {\si} {\sigma}
\newcommand {\Si} {\Sigma}
\newcommand {\de} {\delta}
\newcommand {\prtl} {\partial}
\newcommand {\fr} {\displaystyle\frac}
\newcommand {\wt} {\widetilde}
\newcommand {\be} {\begin{equation}}
\newcommand {\ee} {\end{equation}}
\newcommand {\ba} {\begin{array}}
\newcommand {\ea} {\end{array}}
\newcommand {\bp} {\begin{picture}}
\newcommand {\ep} {\end{picture}}
\newcommand {\bc} {\begin{center}}
\newcommand {\ec} {\end{center}}
\newcommand {\bt} {\begin{tabular}}
\newcommand {\et} {\end{tabular}}
\newcommand {\lf} {\left}
\newcommand {\rg} {\right}

\newcommand {\cF} {{\cal F}}

\newcommand {\cR} {{\cal R}}
\newcommand {\cS} {{\cal S}}

\newcommand {\ses} {\medskip}

\newcommand {\e} {\mathop{\rm e}\nolimits}

\newcommand {\arccosh} {\mathop{\rm arccosh}\nolimits}

\newcommand {\g}  {\stackrel{g\to -g}{\Longleftrightarrow}}

\newcommand {\nin} {\noindent}
\newcommand {\set}[1] {\mathbb{#1}}

\newcommand {\De} {\Delta}

\usepackage{psfrag}

\usepackage{epsfig, latexsym,graphicx}
 \usepackage{amsmath}

\def\2#1#2#3{{#1}_{#2}\hspace{0pt}^{#3}}
\def\3#1#2#3#4{{#1}_{#2}\hspace{0pt}^{#3}\hspace{0pt}_{#4}}
\newcounter{sctn}
\def\sec#1.#2\par{\setcounter{sctn}{#1}\setcounter{equation}{0}
                  \noindent{\bf\boldmath#1.#2}\bigskip\par}

\begin {document}

\begin {titlepage}


 \ses\ses\ses\ses

\vspace{0.1in}

\begin{center}

{\Large \bf Finsleroid--Finsler  Space with Berwald and }
\\
\ses\ses
{\Large \bf Landsberg Conditions}

\end{center}

\vspace{0.3in}

\begin{center}

\vspace{.15in} {\large G.S. Asanov\\} \vspace{.25in}
{\it Division of Theoretical Physics, Moscow State University\\
119992 Moscow, Russia\\
{\rm (}e-mail: asanov@newmail.ru{\rm )}} \vspace{.05in}

\end{center}

\begin{abstract}

We formulate the notion of the  Finsleroid--Finsler space, including the positive--definite as well as indefinite
cases.
 The associated concepts of angle, scalar product, and the distance function
 are elucidated.
If the Finsleroid--Finsler space is of Landsberg type, then  the
Finsleroid charge is a constant. The Finsleroid--Finsler space
proves to be
 a Berwald space if and only if the Finsleroid--axis 1-form is parallel
with respect to the associated Riemannian metric
 and, simultaneously, the Finsleroid charge is a constant.
The necessary and sufficient conditions for the Finsleroid--Finsler space to be of the Landsberg type
are found, which are explicit and simple.
The structure of the associated curvature tensors has been elucidated.

 \ses\ses

\nin

\end{abstract}

\end{titlepage}

\vskip 1cm

\setcounter{sctn}{1}
\setcounter{equation}{0}

 \bc
{\large 1. Introduction and synopsis of  conclusions} \ec

\bigskip

The ideas and methods of the Riemannian geometry admit ingenious   successful metric  extensions
 by using
various  Finslerian  metric functions
(see [1-12]).
The Finsleroid--type metric function $K$,
which properties   have been  systematically exposed
in the previous work [12--14],
is interesting because of entailing numerous remarkable implications,
 including the positive--definiteness,
the constant curvature of the indicatrix,
the special algebraic form of the Cartan tensor, the occurrence of the angle
and scalar product, the particular conformal properties,
the distance  as well as the explicit solutions to the geodesic equations in the respective Finsleroid--Minkowski
space, {\it  etc}.
The Finsleroid is the unit ball defined by the function $K$ under treating  $K$  as a Minkowskian norm.
The respective Finslerian indicatrix  is the boundary of the Finsleroid, which is
strongly   convex and rotund.

Below we
 introduce an appropriate extension  to work on curved manifolds,
such that the Finsleroid--Minkowski structure gets placed on each tangent space.
Namely, given a Riemannian space under assumption that the space admits a 1-form,
 we  use the involved Riemannian metric ${\cal S}$
 and the 1-form, introducing
also the scalar  $g$
which plays the role of  the  Finsleroid charge,
to construct on the Riemannian space
the Finsleroid--Finsler metric function $K$ such that
in each associated tangent space  the vector entering the 1-form  plays the role of the axis of the Finsleroid.
We call  the result the
{\it Finsleroid--Finsler Space},
the {\it $\cF\cF^{PD}_g $--space} for short.
The space thus appeared is actually some structure over a Riemannian space.
 Whenever $g=0$,
 the Finsleroid is the Euclidean unit ball and whence
the  $\cF\cF^{PD}_g $--space
 is a Riemannian space.

Clarifying and studying the entailed Finslerian geometric structure,
 the connection and curvature included, seems to be an urgent task.
First of all, it is of interesting to clarify the Berwald and Landsberg types among the metric functions
of the   $\cF\cF^{PD}_g $--space, for the types play an important role in the modern Finsler geometry
(see [7-9,15]) and sound intriguing to apply in  developments of physical applications.
Accordingly, we start with raising up the following questions.

\ses\ses

QUESTION 1.
{\it Do there exist non--trivial conditions under which
the Finsleroid--Finsler space
is  a Berwald space}?

\ses\ses

QUESTION 2.
{\it Can the Finsleroid--Finsler space  be
of Landsberg type and not of Berwald type}?

\ses\ses


Let
$\{A_{ijk}\}$ denote the Cartan tensor
and
$
A_{i} :=
g^{jk}A_{ijk}
$
be the contraction thereof by means of the Finslerian metric tensor.
Among the conditions
\be
A_{jmn|i}=0,
\ee
\be
A_{j|i}=0,
\ee
\be
\dot A_{jmn}=0,
\ee
\be
\dot A_{j}=0
\ee
\ses
the last one is obviously weakest.
We shall comply with the notation adopted in the book [7], so that
the  bar means the $h$-covariant derivative,  the dot over $A$ stands for the action of the operator
$|il^i$, where $l^i=y^i/K$, and the indices $i,j,\dots$ refer to natural local coordinates
(denoted by $x^i$); the notation $y$ is used for tangent vectors;
\be
\dot X=X_{|i}l^i
\ee
for any tensor $X.$
For purposes of calculations performed below, the smoothness of the class $C^3$ is sufficient for the  associated
Riemannian space
$\cR_N$  as well as  for the input 1-form
$b=b_i(x)y^i$  and the scalar $g(x)$.
The upperscript ``PD" emphasizes the Positive--Definiteness (see the determinant value (2.25) in Section 2)
of
the  $\cF\cF^{PD}_g $--space under study.

A Finsler space is said to be a {\it Berwald space} if the condition (1.1) holds.
A  {\it Landsberg space}  under the Finslerian consideration is  characterized by the  condition (1.3).

\ses\ses

INITIAL OBSERVATIONS.
{\it If the weakly Landsberg condition \rm(1.4)  \it holds, then the Finsleroid charge is independent of points $x$
of the background manifold:
\be
g={\rm constant}.
\ee
The weakly Landsberg condition \rm(1.4)   \it entails the Landsberg condition \rm(1.3), \it and
the weakly Berwald condition \rm(1.2) \it  entails the Berwald condition \rm(1.1).

\ses

 A mere glance at the algebraic structure (2.27)--(2.28) of the   $\cF\cF^{PD}_g $--space
 Cartan tensor
is sufficient to recognize the validity of these assertions.

\ses

Below, we shall confine the treatment to the dimensions
\be
N\ge 3,
\ee
keeping in mind that
the two--dimensional  $\cF\cF^{PD}_g $--space is of  simple structure (we shall present  required notes in the last section
Conclusions).


In our consideration, an important role will be played by the {\it lengthened  angular metric tensor}
\be
{\cal H}_{ij}~:=h_{ij}-\fr{A_iA_j}{A_nA^n}
\ee
which fulfills the identities
\be
{\cal H}_{ij}y^j=0,\qquad {\cal H}_{ij}A^j=0,\qquad{\cal H}_{ij}b^j=0;
\ee
$h_{ij}$ is the traditional Finslerian angular metric tensor (see (2.24)).

We set forth  the following  theorems.

\ses\ses

 {\large Theorem 1.} {\it  Whenever the Finsleroid charge is a constant, the condition} (1.6),
{\it  the representation
\be
{\dot A}_i=
\fr {Ng}{2q}{\cal H}_i{}^mP_m
\ee
takes place,
where
\be
P_m=y^j\nabla_jb_m
+\fr12gqb^j\nabla_jb_m.
\ee

}

\ses\ses

Here,
$q$ is  the quantity (2.3)
and
 $\nabla$ stands for the covariant derivative with respect to the underlined Riemannian metric
 ${\cal S}$.

 Whenever the Finsleroid charge is a constant, all the three identities
\be
\dot A_jy^j=0,\qquad\dot A_jA^j=0,\qquad \dot  A_j b^j=0
\ee
are fulfilled. Indeed, the first identity is universal for any Finsler space.
In the   $\cF\cF^{PD}_g$--space, the contraction $A^hA_h$ is a function of the Finsleroid charge $g$
(see (2.28)), so that whenever
$g=const$
we must have $\dot A_jA^j=0$.
The  third identity is
a direct implications of the first and second ones because  $A_i$
are linear combinations of $y_i$ and $b_i$ (see (A.7) in Appendix A).

Comparing (1.12) with (1.9) gives a handy clue to search for the object ${\dot A}_i$
of  the $\cF\cF^{PD}_g$--space to be a contraction of
the tensor ${\cal H}^k{}_i$ by a vector. The straightforward and attentive (and rather lengthy)
calculations
result in (1.10) and (1.11));
the vanishing
\be
b^j\nabla_ib_j=0
\ee
(the vector $b_i$ is of the unit Riemannian length according to the formula (2.5) of the next section)
nullifies many terms appearing while calculating.
In particular, applying (1.13) to (1.11) results in
\be
b^iP_i=0.
\ee


Let us turn to the Berwald condition.

\ses\ses

 {\large Theorem 2.} {\it   The   $\cF\cF^{PD}_g$--space is a Berwald space if and only if
 in addition to the constancy} (1.6) {\it of the  Finsleroid charge
 the  Finsleroid--axis field $b_i(x)$ is parallel},
 $  \nabla_ib_j=0  $,
{\it with respect to the input Riemannian metric} ${\cal S}$.
}

\ses
\ses

 Under the conditions set forth in this theorem,
 the $hh$-connection and the $hh$-curvature tensor are  the Riemannian connection and the Riemannian curvature tensor
  produced by the input Riemannian metric ${\cal S}$,
and  the $hv$-curvature tensor $P$  vanishes identically;
  at the same time, the basic Finslerian metric tensor remains
 being non--Riemannian (unless $g=0$).

Theorem 2 is a particular case of the following.


\ses\ses

 {\large Theorem 3.} {\it   The   $\cF\cF^{PD}_g$--space is a Landsberg space if and only if
 in addition to the constancy} (1.6) {\it of the  Finsleroid charge
 the  Finsleroid--axis {\rm 1}-form  $b=b_i(x)y^i$
 is closed
\be
\nabla_jb_i-\nabla_ib_j=0
\ee
and
obeying the condition
 \be
\nabla_ib_j=k(a_{ij}-b_ib_j)
 \ee
with a scalar $k=k(x)$.

}

\ses
\ses

In such a space, it proves  possible to obtain a simple explicit representation for each
distinguished Finslerian tensor. In particular, we have the following.

\ses
\ses

 {\large Theorem 4.} {\it  In any Landsberg case of
the $\cF\cF^{PD}_g $--space, the tensor
\be
P_{ijkl}~:=-A_{ijk|l}
\ee
 is of the special algebraic form
\be
P_{ijkl}=
-\fr{gB}{2qK}
\bigl(
{\cal H}_{ij}{\cal H}_{kl}+{\cal H}_{ik}{\cal H}_{jl}+{\cal H}_{jk}{\cal H}_{il}
\bigr)k,
\ee
 and therefore,
}
\be
A_{i|j}=\fr{NgB}{2qK}
{\cal H}_{ij}k.
\ee

\ses
\ses

In (1.16), (1.18),  and (1.19), the case $k=0$  corresponds to the Berwlad type.

In the above formulae,   $K$ is the Finsleroid metric function
given by (2.14)--(2.18) and $B$ is the characteristic quadratic form (2.11).
 The formula (1.18) is remarkable in that
the tensor (1.17) is lucidly constructed in terms of the  lengthened  angular metric tensor
${\cal H}_{ij}$
defined by  (1.8).
Complete  representation for the $hv$-curvature tensor $P_k{}^i{}_{mn}$
can be found in the end of Appendix A below.


The  tensor $P_{ijkl}$ given by
(1.18) is {\it totally symmetric} in all four of its indices and, owing to (1.9), is entailing the identities
\be
P_{ijkl}y^l = 0, \qquad P_{ijkl}A^l = 0, \qquad P_{ijkl}b^l = 0,
\ee
which in turn entail
\be
A_{i|j}y^j=0, \qquad A_{i|j}A^j=0, \qquad A_{i|j}b^j=0.
\ee

The sufficiency in Theorem 3 is easy to verify by direct (two-night-long!) computations.
As regards the necessity,
the representation (1.10)--(1.11) given by Theorem 1 proves to be not only elegant but also handy in permitting
 the equation
${\dot A}_i=0$
 to be resolved,
 using a rather simple   Lemma proven in Section 3.
The shortest way to arrive at the representation (1.18) stated by Theorem 4 is to apply the Finslerian formula
(3.1) of Section 3.

The Landsberg--characteristic condition (1.16) entails the equality
\be
y^iy^j\nabla_jb_i=kq^2
\ee
and also the identity
\be
b^j\nabla_jb_i=0
\ee
(supplementary to (1.13)). If we consider the {\it $b$-lines}
which comprise the congruence tangent to the vector field $b^i(x)$, noting also
the Riemannian unity (see (2.5) below in Section 2) for the length of the vectors $b^i$,
we are justified in concluding from (1.23) that the following proposition is a truth.

\ses \ses

{\large Proposition 1.} {\it  In any  Landsberg--case of
the $\cF\cF^{PD}_g $--space,
 the  $b$-lines are  geodesics of the associated Riemannian space}.

 \ses\ses


When the Landsberg case takes place for a $\cF\cF^{PD}_g $--space,
the computation of the key contraction
$\ga^i{}_{nm}y^ny^m$,  with $
\ga^i{}_{nm}
$
being the associated Finslerian Christoffel symbols,
leads us  to the following astonishingly simple result:
\be
\ga^i{}_{nm}y^ny^m=gqk(y^i-bb^i)
+a^i{}_{nm}y^ny^m,
\ee
where
$a^i{}_{km}$
are
the Christoffel symbols constructed from the input Riemannian metric ${\cal S}$.

Since the representation  (1.24) entails the nullification
\be
b_i(\ga^i{}_{nm}y^ny^m-a^i{}_{nm}y^ny^m)=0
\ee
(notice that $b^ib_i=1$ in accordance with the formulas (2.5) and (2.8) of the next section),
we are justified in   stating  the following remarkable assertion.

\ses \ses

{\large Proposition 2.} {\it  In any Landsberg case of
the $\cF\cF^{PD}_g $--space,
the covariant derivatives of the vector field $b_i(x)$ with respect to the Finsler connection
and with respect to the associated Riemannian connection are equal to one another}:
\be
b_{i|j}=\nabla_jb_i.
\ee

\ses\ses

The commutator
 \be
(\nabla_m\nabla_n-\nabla_n\nabla_m)b_k=\wt k_m(a_{nk}-b_nb_k)-\wt k_n(a_{mk}-b_mb_k)
\ee
(insert the Landsberg--characteristic condition (1.16)  in the left--hand part)
is valid with the vector
\be
\wt k_n=\D k{x^n}+k^2b_n.
\ee

\ses\ses

Finally, the formulas obtained in the end of Appendix A enable us to formulate the following.

\ses\ses

{\large Proposition 3.} {\it  In any Landsberg case of
the $\cF\cF^{PD}_g $--space,
the tensor
\be
\rho_{ij}~:=\fr12(R_i{}^m{}_{mj}+R^m{}_{ijm})-\fr12g_{ij}R^{mn}{}_{nm}
\ee
is covariantly conserved:
\be
\rho^i{}_{j|i} \equiv 0.
\ee
}


Lengthy straightforward calculations lead to the following assertion.

\ses\ses

{\large Proposition 4.} {\it  Let us be given a $\cF\cF^{PD}_g $--space with constant Finsleroid charge.
Then the tensor ${\dot A}_{ikl}$ is of the following special algebraic structure:
\be
{\dot A}_{ikl}=\fr g{2q}\Bigl(c_i{\cal H}_{kl}+c_k{\cal H}_{li}+c_l{\cal H}_{ik}\Bigr)
\ee
 with
\be
c_i=\Bigl(y^j\nabla_jb_i
+\fr12gqb^j\nabla_jb_i\Bigr)
-\fr1{q^2}\Bigl(y^jy^h\nabla_jb_h+\fr12gqy^hb^j\nabla_jb_h\Bigr)v_i,
\ee
where $v_i=a_{in}y^n-bb_i$.
}

\ses\ses

This just entails the representation
\be
\dot A_i=\fr{gN}{2q}c_i
\ee
which is equivalent to (1.10). The vector
(1.32) fulfills the identities
\be
c_iy^i=0, \qquad c_ib^i=0, \qquad c_iA^i=0
\ee
(notice that $v_ib^i=0$).


Below, in  Section 2 we formulate rigorously  the basic concepts
which underline the notion of the  $\cF\cF^{PD}_g $--space.
The Finsleroid--Finsler metric function, $K$, given by the explicit formulas (2.14)--(2.18),
is defined by the triple: a Riemannian metric $S$, a 1-form $b$, and a scalar $g(x)$.
We attribute to the vector $b^i(x)$ of the 1-form the geometrical meaning of the direction of the
 {\it axis} of the Finsleroid supported by a point $x\in M$, and to
the scalar $g(x)$  the meaning of ``the geometric charge", particularly
``the Finsleroid charge". Would the charge $g(x)$ be zero, the function $K$ reduce to the Riemannian $S$.
The special and lucid algebraic form (2.27)--(2.28) of the associated Cartan tensor
gives rise to many essential simplifications, including  the constant curvature of the
Finsleroid indicatrix.
The explicit formula (2.34) is proposed for the angle that equals  the  length of a piece of the Finslerian
unit circle on the Finsleroid indicatrix (the meaning of such an angle complies
 totally with the ordinary meaning
of the Euclidean ``induced angle'' on the tangent spheres in the Riemannian geometry).
Having obtained the angle, we obtain the scalar product, too.

In Section 3 we place the arguments which verify the validity of the above theorems.

In Section 4
the structure of the associated curvature tensor is elucidated in various aspects.

In Section 5
we open up the remarkable phenomenon that the Landsberg--type   $\cF\cF^{PD}_g $--space  is underlined by
a warped product metric of the associated Riemannian space with respect to the $b$-geodesic
coordinates.

In Section 6 we indicate a convenient possibility to re--formulate the theory in  the indefinite case
of the time--space type.
The resultant space will be referred to as the  $\cF\cF^{SR}_g $--space.
Although the  above formulas and the analysis presented in Sections 3--5 are all referred to
 the positive--definite case,
a simple additional consideration (omitted in the present paper) shows that
all the principal representations and conclusions
 can straightforwardly be re--addressed to the indefinite pseudo--Finsleroid space exposed in Section 6:
to do this, in most cases  it is sufficient merely to change in formulas the sign  of the Finsleroid charge $g$.
Whence all Theorems 1-7 and Propositions 1--5 remain valid under  transition from the positive--definite
Finsleroid--Finsler  $\cF\cF^{PD}_g $--space to the indefinite
pseudo--Finsleroid--Finsler   $\cF\cF^{SR}_g $--space.

 The basic geometric ideas involved
 and possible applied  potentialities of the spaces
 proposed are emphasized in Conclusions.

  Appendix A summarizes the important representations  which  support the calculations involved.


\setcounter{sctn}{2}
\setcounter{equation}{0}

\bc
{\large 2. Positive--Definite Case}
\ec

\ses\ses

Let $M$ be an $N$-dimensional
differentiable  manifold.
Suppose we are given on $M$ a Riemannian metric ${\cal S}=S(x,y)$, where $x\in M$ denotes  points
and $y\in T_xM$ means tangent vectors. Denote by
$\cR_N=(M,{\cal S})$
the obtained $N$-dimensional Riemannian space.

Let us assume that the manifold $M$ admits a non--vanishing 1-form, to be denoted as $\beta=\beta(x,y)$, and
introduce the associated normalized 1-form $b=b(x,y)$
 such that the Riemannian length of the involved vector   be equal to 1.
With respect to  natural local coordinates in the space
$\cR_N$
we have the local representations
\be
 b=b_i(x)y^i,
\ee
\ses
\be
S= \sqrt{a_{ij}(x)y^iy^j},
\ee
\ses
\be
q=\sqrt{r_{ij}(x)y^iy^j}
\ee
with
\be
r_{ij}(x)=a_{ij}(x)-b_i(x)b_j(x)
\ee
and
\be
a^{ij}b_ib_j=1.
\ee
The decomposition
\be S^2=b^2+q^2
\ee
introduces a Riemannian metric $q$ on the $(N-1)$--dimensional Riemannian space $\cR_{N-1}\subset \cR_N$
 given rise to by the orthogonality relative to the 1-form  $b$, for
from (2.4) and (2.5) it follows that
\be
b^ir_{ij}=0,
\ee
where
\be
b^i~:=a^{ik}b_k.
\ee

 Finally, we introduce on $M$
a scalar field $g=g(x)$
 subject to ranging
\be
-2<g(x)<2,
\ee
and apply  the convenient notation
\be
h=\sqrt{1-\fr14g^2}, \qquad
G=g/h.
\ee
\ses


The {\it  characteristic
quadratic form}
\be
B(x,y) :=b^2+gqb+q^2
\equiv\fr12\Bigl[(b+g_+q)^2+(b+g_-q)^2\Bigr]>0
\ee
where $ g_+=\fr12g+h$ and $ g_-=\fr12g-h$,
is of the negative discriminant
\be
D_{\{B\}}=-4h^2<0
\ee
and, therefore, is positively definite.
In the limit $g\to 0$,
the definition (2.11) degenerates to the
 quadratic form (2.6) of the input Riemannian metric tensor:
\be
B|_{_{g=0}}=S^2.
\ee

In terms of these concepts, we extend the notion of the Finsleroid metric function $K$ proposed
early  in the framework of the Minkowski space (see [12--14]) and introduce the following
definition adaptable to consideration on manifolds.

\ses\ses

 {\large Definition}. The scalar function $K(x,y)$ given by the formulae
\be
K(x,y)=
\sqrt{B(x,y)}\,J(x,y)
\ee
and
\be
J(x,y)=\e^{\frac12G\Phi(x,y)},
\ee
where
\be
\Phi(x,y)=
\fr{\pi}2+\arctan \fr G2-\arctan\Bigl(\fr{L(x,y)}{hb}\Bigr),
\qquad  {\rm if}  \quad b\ge 0,
\ee
and
\be
 \Phi(x,y)= -\fr{\pi}2+\arctan \fr
G2-\arctan\Bigl(\fr{L(x,y)}{hb}\Bigr), \qquad  {\rm if}  \quad
b\le 0,
\ee
 with
 \be
 L(x,y) =q+\fr g2b,
\ee
\ses\\
is called
the {\it  Finsleroid--Finsler  metric function}.

\ses\ses

The positive (not absolute) homogeneity holds fine: $K(x,\la y)=\la K(x,y)$ for all $\la >0$.

Sometimes it is convenient to use also the function
\be
A(x,y)=b+\fr g2q.
\ee
We have
\be
L^2+h^2b^2=B, \qquad A^2+h^2q^2=B.
\ee


\ses\ses

 {\large  Definition}.  The arisen  space
\be
\cF\cF^{PD}_g :=\{\cR_{N};\,b(x,y);\,g(x);\,K(x,y)\}
\ee
is called the
 {\it Finsleroid--Finsler space}.

\ses\ses

 {\large  Definition}. The space $\cR_N$ entering the above definition is called the {\it associated Riemannian space}.

\ses\ses

{\large Definition}.\, Within  any point $T_xM$, the Finsleroid--metric function $K(x,y)$  produces the {\it Finsleroid}
 \be
 \cF^{PD}_{g\,\{x\}}:=\{y\in  \cF^{PD}_{g\,\{x\}}: y\in T_xM , K(x,y)\le 1\}.
  \ee

\ses \ses

{\large Definition}.\, The {\it Finsleroid Indicatrix}
 $
I^{PD}_{g\,\{x\}}\in T_xM$ is the boundary of the Finsleroid:
 \be
I^{PD}_{g\,\{x\}} :=\{y\in I^{PD}_{g\,\{x\}}: y\in T_xM, K(x,y)=1\}.
  \ee

\ses \ses

 Since at $g=0$ the  $\cF\cF^{PD}_g$--space is
Riemannian, then the body $  \cF^{PD}_{g=0\,\{x\}}$ is a unit ball and $
I^{PD}_{g=0\,\{x\}}$ is a unit sphere.

\ses
\ses

 {\large  Definition}. The scalar $g(x)$ is called
the {\it Finsleroid charge}.
The 1-form $b$ is called the  {\it Finsleroid--axis}  1-{\it form}.

\ses\ses

Under these conditions,
 we can explicitly calculate from the function $K$ the  distinguished Finslerian tensors,
 and first of all
the covariant tangent vector $\hat y=\{y_i\}$,
the  Finslerian metric tensor $\{g_{ij}\}$
together with the contravariant tensor $\{g^{ij}\}$ defined by the reciprocity conditions
$g_{ij}g^{jk}=\de^k_i$, and the  angular metric tensor
$\{h_{ij}\}$, by making  use of the following conventional  Finslerian  rules in succession:
\be
y_i :=\fr12\D{K^2}{y^i}, \qquad
g_{ij} :
=
\fr12\,
\fr{\prtl^2K^2}{\prtl y^i\prtl y^j}
=\fr{\prtl y_i}{\prtl y^j}, \qquad
h_{ij} := g_{ij}-y_iy_j\fr1{K^2}.
\ee
Calculations show that the determinant of the associated Finsleroid metric tensor  is everywhere positive:
\be
\det(g_{ij})=\bigl(\fr{K^2}B\bigr)^N\det(a_{ij})>0.
\ee


After that, we can elucidate the algebraic structure of the associated  Cartan tensor
\be
A_{ijk} := \fr K2\D{g_{ij}}{y^k},
\ee
which leads to
the following simple and remarkable result:
\it The Cartan tensor associated with the  Finsleroid--Finsler  metric function $K$
 is of the following special algebraic form:
\be
A_{ijk}=\fr1N\lf(h_{ij}A_k+h_{ik}A_j+h_{jk}A_i-\fr1{A_hA^h}A_iA_jA_k\rg)
\ee
with
\be
A_hA^h=\fr{N^2}{4}g^2,
\ee
\rm
which is in fact a mere adaptation of the
 Cartan tensor  structure known in the Finsleroid--Minkowski approach developed in  [12--14].
Elucidating the  respective  tensor
\be
\hat R_i{}^j{}_{mn} := \fr1{K^2}(\3Ahjm\3Aihn-\3Ahjn\3Aihm)
\ee
describing the curvature of the indicatrix results, upon using (2.27),  in the simple representation
$$
\hat R_{ijmn}=-\fr{A_hA^h}{N^2}(h_{im}h_{jn}-h_{in}h_{jm}).
$$
Inserting here (2.28), we are led to the following remarkable assertion.
\it The indicatrix curvature tensor \rm(2.29) \it of the space
$\cF\cF^{PD}_g$ is of the special structure such that
\be
\hat R_{ijmn}=S^*(h_{im}h_{jn}-h_{in}h_{jm})/K^2
\ee
with \rm
\be
S^*=-\fr14g^2.
\ee

Recalling the known formula $ \cR=1+S^* $ for the indicatrix
curvature (see Section 1.2 in [1]), from (2.29) and (2.30)  we conclude that
\be
 \cR_{\text{Finsleroid  Indicatrix} }=h^2\equiv 1-\fr14g^2,
 \ee
  so that $$ 0
< \cR_{\text{Finsleroid Indicatrix }} \le 1 $$ and
$$
\cR_{\text{Finsleroid Indicatrix} }\stackrel{g\to 0}{\Longrightarrow}
\cR_{\text{Euclidean Sphere}}=1.
$$
 Geometrically, the fact that the quantity ~(2.32)
is independent of  vectors~$y$ means:
 \it The Finsleroid indicatrix
 is a space of constant  positive curvature\rm.

\ses
\ses


Given  any two nonzero tangent vectors
  $y_1,y_2\in T_xM$ of a fixed tangent space, we can, by following the previous work [12--14],
obtain the {\it $\cF\cF_g^{PD}$--scalar product}
 \be
<y_1,y_2>_{\{x\}}~:=K(x,y_1)K(x,y_2) \cos\Bigl( \al_{\{x\}}(y_1,y_2)\Bigr)
 \ee
 and the $\cF\cF_g^{PD}$--{\it angle}
 \be
 \al_{\{x\}}(y_1,y_2)~: = \fr1h\arccos \fr{ A(x,y_1)A(x,y_2)+h^2<y_1,y_2>_{\{x\}}^{\{r\}} }
 {\sqrt{B(x,y_1)}\,\sqrt{B(x,y_2)} },
 \ee
where $<y_1,y_2>_{\{x\}}^{\{r\}}=r_{ij}(x)y_1^iy_2^j$.
The $\cF\cF_g^{PD}$--{\it distance} $\De s$ between ends of the vectors in the tangent space
is given by the formula
\be
 (\Delta s)^2 = (K(x,y_1))^2 + (K(x,y_2))^2 - 2 K(x,y_1)K(x,y_2) \cos\Bigl( \al_{\{x\}}(y_1,y_2)\Bigr),
 \ee
which
extends the ordinary Euclidean cosine theorem.
At equal vectors, the reduction
\be
 <y,y>_{\{x\}}=(K(x,y))^2
  \ee
 takes place, that is, the two-vector scalar product (2.33)
reduces  exactly  to the squared Finsleroid--Finsler metric function.

These scalar product and angle obviously exhibit the symmetry:
\be
(y_1,y_2)_{\{x\}}=(y_2,y_1)_{\{x\}}, \qquad
\al_{\{x\}}(y_1,y_2)= \al_{\{x\}}(y_2,y_1).
\ee
They are   entirely intermediary, supporting  by a point $x\in M$ of the base manifold $M$
(in just the  same sense as in the Riemannian geometry)
and being  independent of any vector element of support.
The angle (2.34) furnishes actually the ordinary meaning of the arc--length $s$ cut off by a geodesic piece
on the indicatrix:
\be
\al_{\{x\}}(y_1,y_2)=
s_{\{\rm Finsleroid~ Indicatrix\}}(x;y_1,y_2).
\ee
Namely, considering the  indicatrix together with the geodesic curves thereon,
we can calculate
 the length of the geodesic piece  which joins the ends of the
 unit vectors
$l_1=y_1/K(x,y_1)$ and
$l_1=y_2/K(x,y_2)$,
thereby obtaining the right--hand part in (2.34).
Taken (2.38) to be  the definition for the  angle,
 the explicit representation (2.34) can be verified by direct calculation
of the arc length over the Finsleroid indicatrix;
the same angle is derivable by postulating the cosine theorem (2.35)
(see the previous work
[12--14]).
The formula (2.34) assigns the angle  to be  a length of a piece of the Finslerian
unit circle on the Finsleroid indicatrix.
Whenever the Finsleroid charge is zero, $g=0$, the Finslerian angle (2.34) reduces to the ordinary Riemannian angle
$ \arccos\Bigl[\bigl(b(x,y_1)b(x,y_2)+<y_1,y_2>_{\{x\}}^{\{r\}}\bigr)\,\Bigl(S(x,y_1)S(x,y_2)\Bigr)^{-1}\Bigr]
=
 \arccos\Bigl[a_{ij}(x)y^i_1y^j_2\,\Bigr(S(x,y_1)S(x,y_2)\Bigl)^{-1}\Bigr].$


In particular, if we consider the angle
$ \al_{\{x\}}(y)$ formed by a vector $y\in T_xM$ with the Finsleroid axis in a fixed $T_xM$,
from (2.33) we get the respective value to be
 \be
 \al_{\{x\}}(y)~: = \fr1h\arccos \fr{ A(x,y)}
 {\sqrt{B(x,y)}}.
 \ee

\setcounter{sctn}{3}
\setcounter{equation}{0}

\bc
{\large 3. Search for Berwald and Landsberg types of $\cF\cF^{PD}_g $--space}
\ec

\ses\ses

In any Finsler space subjected to the Landsberg--type condition
the equality
\be
A_{ijk|l}=-\Bigl(F\D{A_{ijk}}{y^l}\Bigr)_{|s}l^s
\ee
(see
(12.4.1) on p. 326 of the book [7])
holds.
In the $\cF\cF^{PD}_g $--space under study, the representations (A.23) and (A.24)  may be
used to  just conclude from (3.1) that
\be
{\cal A}_{ijk|l}=
\fr{g}{2}\bigl(
{\cal H}_{ij}{\cal H}_{kl}+{\cal H}_{ik}{\cal H}_{jl}+{\cal H}_{jk}{\cal H}_{il}
\bigr)\dot {\la},\qquad \la=\fr bq,
\ee
in terms of the lengthened  angular metric tensor
$
{\cal H}_{ij}
$
(defined by  (1.8)).

The representation (3.2) is beautiful but only suggestive, assigning  a necessary condition.
To gain the full sufficient condition for
the $\cF\cF^{PD}_g $--space  to be of a Landsberg type,
we are to straightforwardly calculate  the object
$ \dot A_i $, taking into account the constancy (1.6) of the Finsleroid charge and the vanishing
$
b^j\nabla_ib_j=0
$
(see (1.13)). On so doing, we arrive at
merely
\be
{\dot A}_i=
\fr{1}{q^2K}
\Bigl[(b+gq)A_i
-\fr N2
gq
l_i
\Bigr]
y^kP_k
+\fr N2\fr{g}q
P_i
\ee
with $P_i$ being given by (1.11).
This representation (3.3) can readily be rewritten in the form (1.10).


It is easy to prove the following.

\ses\ses

{\large Lemma.} {\it In  the  $\cF\cF^{PD}_g$--space  with}
  $ g=const\ne 0$,
 \be
 \dot A_i=0 \quad \Longleftrightarrow \quad \nabla_jb_i=k(x)(a_{ij}-b_ib_j).
 \ee

\ses\ses

\ses\ses

Indeed, the right--hand side of the obtained representation (3.3) for the object
$\dot A_i$ is such that
the vanishing $ \dot A_i=0$
takes place if and only if the vector $P_i$ is a linear combination of $u_i=a_{ij}y^j$ and $b_i$:
\be
P_i=k_1u_i+k_2b_i
\ee
(the tensor ${\cal H}_{ik}$ is of the rank $N-2$ and at each point $x\in M$
has the isotropic plane spanned by
the vectors $u_i=a_{ij}y^j$ and $b_i$ because of the identities (1.9)).
To agree with the identity $b^iP_i=0$ (see (1.14)), we must put  $k_2=-bk_1$, obtaining
\be
y^j\nabla_jb_i
+\fr12gqb^j\nabla_jb_i=k_1(u_i-bb_i)\equiv k_1r_{ij}y^j.
\ee
If we apply here the operator $b^k\partial/\partial y^k$
 and take into account  the vanishing $r_{ij}b^j=0$
 (see (2.7)), together with
 $b^k\partial q/\partial y^k=0$ (which is tantamount to (2.7)), we may conclude that
$
b^j\nabla_jb_i=0,
$
whence in (3.6) the term $
y^j\nabla_jb_i
$
must be proportional to $u_i-bb_i$.
Therefore, the examined condition  $ \dot A_i=0$ reduces to
$
y^j\nabla_jb_i
=k(u_i-bb_i),
$
where $k$ is a factor. Since here the vectors
$
y^j\nabla_jb_i$
and
$
u_i-bb_i=(a_{ij}-b_ib_j)y^j$
are linear functions of the set $\{y^j\}$, the factor $k$ may not depend on $y^j$.
Therefore,
$
y^j\nabla_jb_i
=k(x)(u_i-bb_i).
$
Differentiating the last result with respect to $y^j$ just yields
$
\nabla_jb_i
=k(a_{ij}-b_ib_j)$ with
$ k=k(x).
$
 Lemma is valid.

\ses\ses

Therefore,  Theorems 3 and 4 are valid.


In processing calculations, it is frequently necessary to keep in mind that the formulas (1.24)--(1.26)
entail the following.

\ses \ses

{\large Proposition 5.} {\it  In any Landsberg case of
the $\cF\cF^{PD}_g $--space,
the covariant derivative of the input Riemannian metric tensor $a_{mn}$ is such that
$$
a_{mn|k}=
-\fr {gk}q\biggl[
(a_{km}-b_kb_m)(u_n-bb_n)
+
(a_{kn}-b_kb_n)(u_m-bb_m)
+(a_{nm}-b_nb_m)(u_k-bb_k)
$$
\ses
\be
-
\fr 1{q^2}(u_m-bb_m)(u_n-bb_n)(u_k-bb_k)
\biggr].
\ee
For the vector
$u_m=a_{mn}y^n$ we get
\be
u_{m|k}=a_{mn|k}y^n=
-\fr {gk}q\biggl[
q^2(a_{km}-b_kb_m)
+
(u_{k}-bb_k)(u_m-bb_m)
\biggr]
=u_{k|m},
\ee
after which we obtain
\be
b_{|k}=y^n\nabla_kb_n=k(u_{k}-bb_k),
\ee
\ses
\be
(q^2+b^2)_{|k}=a_{mn|k}y^my^n=
-2gq
\nabla_nb,
\ee
\ses
\be
q_{|k}=
-\fr1q(b+gq)\nabla_kb,
\ee
\ses
\be
(bq)_{|k}=\fr1q
[q^2-b(b+gq)]\nabla_kb, \qquad (q/b)_{|k}=-\fr{B}{qb^2}\nabla_kb,  \qquad (b/q)_{|k}=\fr{B}{q^3}
\nabla_kb,
\ee
and
\be
J_{|n}=\fr g{2q}J\nabla_nb, \qquad B_{|n}=-\fr{gB}{q}
\nabla_nb.
\ee
This leads to
\be
K\dot b=kq^2, \qquad  \dot \la=\fr{B}{q^3}\dot b,
\ee
and
\be
S\dot S=-gq\dot b, \quad \dot q=-\fr1q(b+gq)\dot b, \quad \dot{(bq)}=\fr1q [q^2-b(b+gq)]\dot b,
\quad \dot B=-\fr{gB}{q}\dot b.
\ee

}

\ses\ses

It is the implication (3.14) that turns the suggestive representation (3.2) into the conclusive result (1.18).

We have
\be
B_{|m|n}-B_{|n|m}=0\quad \Longleftrightarrow \qquad A_kR^k{}_{mn}=0.
\ee


\setcounter{sctn}{4}
\setcounter{equation}{0}

\bc
{\large 4. Study of Curvature Tensor}
\ec

\ses\ses

When straightforwardly calculating
 the $hh$-curvature tensor $R^i{}_k$ on the basis of the known  definition (see (A.28) in Appendix A below)
with $
\ga^i{}_{nm}y^ny^m$
given by (1.24), we obtain the following.

\ses\ses

{\large Theorem 5.} {\it   In the Landsberg--case
  $\cF\cF^{PD}_g$--space   the following explicit and simple representation is valid:
\be
K^2R^i{}_k=
\fr{g^2k^2}{4}q^2{\cal H}^i{}_k
+gqv^i\wt k_k
-
\fr{g}{2q}(\wt k_ny^n)
(
v^iv_k
+
q^2r^i{}_k
)
+y^na_n{}^i{}_{km}y^m.
\ee

\ses\ses
}

Here,
$v^i=y^i-bb^i$ and $v_m=u_m-bb_m=r_{mn}y^n$
with $u_m=a_{mn}y^n$; $r_{mn}=a_{mn}-b_mb_n$ being the tensor
(2.4);
$r^i{}_n=a^{im}r_{mn}=\de^i{}_n-b^ib_n;$
$a_m{}^i{}_{kn}$ is the Riemannian curvature tensor associated with the input Riemannian metric ${\cal S}$;
$\wt k_n$ is the vector (1.28).
The Berwald case in this theorem is conditioned by $k=0$.

From (4.1) we find the tensor
\be
KR^i{}_{km}=
\fr{g^2k^2}{4}
(
r^i{}_kv_m-
v_kr^i{}_m)
+\fr{ g}{2q}v^i(\wt k_kv_m-\wt k_mv_k)
-\fr12  gq
(r^i{}_k
\wt k_m
-
r^i{}_m\wt k_k
)
+y^na_n{}^i{}_{km}
\ee
(the use of the formula (A.18) of Appendix A is convenient),
and
then the full $hh$-curvature tensor
\ses\\
$$
R_n{}^i{}_{km}=
\fr{g^2k^2}{4}
(
r_{nm}r^i{}_k-r_{nk}r^i{}_m
)
- \fr g{2q}
v_n\Bigl(
{\cal H}^i{}_k
\wt k_m
-
{\cal H}^i{}_m
\wt k_k
\Bigr)
$$
\ses
\be
+\fr{ g}{2q}\Bigl[r^i{}_n(\wt k_kv_m-\wt k_mv_k)+v^i(\wt k_kr_{mn}-\wt k_mr_{kn})\Bigr]
+a_n{}^i{}_{km}
\ee
(apply the rules (A.29) and (A.30) of Appendix A below).
The first term in the right--hand side of (4.3)
is independent of vectors $y^i$.
The equalities
\be
R{}^i{}_{km}=l^nR_n{}^i{}_{km}
\ee
and
\be
R{}^i{}_{k}=l^nR_n{}^i{}_{km}l^m
\ee
hold.


From (4.1) it follows that the {\it Ricci scalar} $Ric~:=R^i{}_i$ is of the value
\be
Ric=
\fr1{K^2}
\Bigl[
\fr{g^2k^2}{4}q^2(N-2)
-\fr12gqN(\wt k_ny^n)
+gq[\wt k_jy^j-b(b^j\wt k_j)]
+y^na_n{}^i{}_{im}y^m
\Bigr];
\ee
from (4.3) it follows that  the tensor
\be
R_{ik}~:=l^nR_{nikm}l^m
\ee
can explicitly be written as
\be
 BR_{ik}=
\Bigl(
\fr{g^2k^2}{4}
-\fr g{2q}(\wt k_my^m)
\Bigr)
(
q^2r_{ik}-v_iv_k
)
+
y^na_{nikm}y^m
\ee
(with
$a_{nikm}=a_{ij}a_n{}^j{}_{km}$),
which can also be written in terms of the lengthened angular metric tensor as follows:
\be
K^2R_{ik}=
\Bigl(
\fr{g^2k^2}{4}
-\fr g{2q}(\wt k_my^m)
\Bigr)
q^2{\cal H}_{ik}
+
y^na_{nikm}y^m\fr {K^2}B.
\ee
The symmetry
\be
R_{ik}=R_{ki}
\ee
and the identity
\be
y^iR_{ik}=0
\ee
are obviously valid.


From (4.3) the   contracted tensor is found:
$$
R_n{}^i{}_{im}=
\fr{g^2k^2}{4}
(N-2)r_{nm}
-\fr12 \fr gqNv_n\wt k_m
$$
\ses
\be
+\fr g{2q}\Bigl[[\wt k_n-b_n(b^j\wt k_j)]v_m
+
[\wt k_m-b_m(b^j\wt k_j)]v_n\Bigr]
+
\fr g{2q}\bigl[(y^j\wt k_j)
-b(b^j\wt k_j)\bigr](r_{nm}-\fr1{q^2}v_nv_m)
+a_n{}^i{}_{im},
\ee
entailing
the following expression for the total contraction of the curvature tensor:
$$
R^{nm}{}_{mn}=g^{nm}
R_n{}^i{}_{im}=
\fr1{K^2}\Biggl[
\fr{g^2k^2}{4}
(N-2)\Bigl((N-1)B+g^2q^2\Bigr)
$$
\ses
\be
-\fr12g^2(N-2)(b+gq)[\wt k_jy^j-b(b^j\wt k_j)]
+\fr12g^2q^2N(\wt k_hb^h)\Biggr]+
a_n{}^i{}_{im}g^{nm}
\ee
with
$$
a_n{}^i{}_{im}g^{nm}=
\fr B{K^2}\Bigl(a_n{}^i{}_{im}a^{nm}-\fr {gb}q(N-1)(\wt k_hb^h)+\fr{2g}q[(N-1)(\wt k_hy^h)-(\wt k_hv^h)]\Bigr)
$$
\ses
\be
+\fr1{K^2}\fr gq(b+gq)y^na_n{}^i{}_{im}y^m.
\ee
We have
\be
R_n{}^i{}_{im}y^ny^m=K^2Ric.
\ee
\ses
The tensor (4.12) is non--symmetric:
\be
R_n{}^i{}_{im}-R_m{}^i{}_{in}=-\fr12\fr gqN(v_n\wt k_m-v_m\wt k_n).
\ee
The contraction
$$
R_n{}^i{}_{im}y^m=
\fr{g^2k^2}{4}
(N-2)v_n
-\fr12 \fr gqNv_n(\wt k_my^m)
$$
\ses
\be
+\fr g{2q}\Bigl[[\wt k_n-b_n(b^j\wt k_j)]q^2
+
[\wt k_my^m-b(b^j\wt k_j)]v_n\Bigr]
+
a_n{}^i{}_{im}y^m
\ee
is obtained.


Let us consider also the {\it Ricci tensor}
\be
Ric_{nm}~:=
\fr12\Dd{(K^2Ric)}{y^n}{y^m},
\ee
\ses
the {\it Ricci--deflection vector}
\be
\Upsilon_n~:=\fr12\D{(K^2Ric)}{y^n}-R_n{}^i{}_{im}y^m
\ee
and
the {\it Ricci--deflection tensor}
\be
\Upsilon_{nm}~:=Ric_{nm}-\fr12(R_n{}^i{}_{im}+R_m{}^i{}_{in})\equiv \Upsilon_{mn}.
\ee
We obtain upon direct calculations the following lucid representations in terms of the lengthened angular metric tensor:
\be
\Upsilon_n=
-\fr14 gqN{\cal H}_n{}^m\wt k_m
\ee
and
\be
\Upsilon_{nm}=
-\fr14\fr gqN(\wt k_jy^j){\cal H}_{nm}.
\ee
The identities
\be
\Upsilon_ny^n=0, \qquad \Upsilon_nA^n=0, \qquad \Upsilon_nb^n=0
\ee
and
\be
\Upsilon_{nm}y^m=0,\qquad\Upsilon_{nm}A^m=0,\qquad\Upsilon_{nm}b^m=0
\ee
are valid.
It is useful to compare (4.21) and (4.22) with (1.10) and (1.19).


To lower the index in the curvature tensor (4.3) we must apply the explicit expression (see (A.2) in Appendix A below)
for the Finsleroid metric tensor.
On so doing, and using the contractions
$$
b_jR_n{}^j{}_{km}=b_ja_n{}^j{}_{km}=-\wt k_kr_{nm}+\wt k_mr_{nk}
$$
(see (1.27)) and
\ses
$$
u_jR_n{}^j{}_{km}=\fr{g^2k^2}{4}
(
r_{nm}v_k-r_{nk}v_m
)
+\fr12gq[\wt k_k(r_{nm}+\fr1{q^2}v_nv_m)-\wt k_m(r_{nk}+\fr1{q^2}v_nv_k)]
+u_ja_n{}^j{}_{km},
$$
we obtain the result
\ses\\
$$
\fr B{K^2}R_{nikm}=
\fr{g^2k^2}{4}
(
r_{nm}r_{ik}-r_{nk}r_{im}
)
-\fr12  \fr gq
v_n\Biggl(
(r_{ik}-\fr1{q^2}v_iv_k)
\wt k_m
-
(r_{im}-\fr1{q^2}v_iv_m)
\wt k_k
\Biggr)
$$
\ses
$$
+\fr{ g}{2q}\Bigl[r_{in}(\wt k_kv_m-\wt k_mv_k)+v_i(\wt k_kr_{mn}-\wt k_mr_{kn})\Bigr]
+a_{nikm}
$$
\ses
$$
+\fr g{qB}
\Biggl[
(gq^3b_i
+S^2v_i)
(-\wt k_kr_{nm}+\wt k_mr_{nk})
$$
\ses\ses
$$
+(S^2b_i-bu_i)
\Bigl(
\fr{g^2k^2}{4}
(
r_{nm}v_k-r_{nk}v_m
)
+\fr12gq[\wt k_k(r_{nm}+\fr1{q^2}v_nv_m)-\wt k_m(r_{nk}+\fr1{q^2}v_nv_k)]
+u_ja_n{}^j{}_{km}
\Bigr)
\Biggr],
$$
which can be  simplified to read
\ses\\
$$
\fr B{K^2}R_{nikm}=
\fr{g^2k^2}{4}
\Bigl[
(
r_{nm}r_{ik}-r_{nk}r_{im}
)
+\fr g{Bq}(S^2b_i-bu_i)
(
r_{nm}v_k-r_{nk}v_m
)
\Bigr]
$$
\ses
$$
+\fr{ g}{2q}\Bigl[(r_{in}-\fr1{q^2}v_iv_n)(\wt k_kv_m-\wt k_mv_k)
-v_i(\wt k_kr_{mn}-\wt k_mr_{kn})\Bigr]
-\fr12  \fr gq v_n(r_{ik} \wt k_m -  r_{im}\wt k_k)
$$
\ses
$$
-\fr {g^2}{2B}
(S^2b_i-bu_i)\Bigl[
\wt k_k(r_{nm}-\fr1{q^2}v_nv_m)-\wt k_m(r_{nk}-\fr1{q^2}v_nv_k)]
\Bigr]
$$
\ses
\be
+\fr g{qB}
(S^2b_i-bu_i)y^ja_{njkm}
+
a_{nikm}.
\ee


Let us use this representation to evaluate the contracted tensor
$g^{nm}R_{nikm}$ (the contravariant components of the metric tensor involved are given by (A.3) in Appendix A below).
We get in succession:
$$
\fr B{K^2}a^{nm}R_{nikm}=
\fr{g^2k^2}{4}(N-2)
\Bigl[
r_{ik}
+\fr g{Bq}(S^2b_i-bu_i)
v_k
\Bigr]
$$
\ses
$$
+\fr{ g}{2q}\Bigl[
-(r_{i}{}^m-\fr1{q^2}v_iv^m)\wt k_mv_k
-v_i(\wt k_k(N-2)-\wt k_mr^m{}_k)\Bigr]
-\fr12  \fr gq
r_{ik}
(\wt k_mv^m)
$$
\ses
$$
-\fr {g^2}{2B}
(S^2b_i-bu_i)\Bigl[
\wt k_k(N-2)-(r_k{}^m-\fr1{q^2}v_kv^m)\wt k_m
\Bigr]
$$
\ses
$$
+
a^{nm}a_{nikm}
+\fr g{qB}
(S^2b_i-bu_i)y^ja_{njkm}a^{nm}
$$
and
\ses\\
$$
\fr B{K^2}y^my^nR_{nikm}=
\fr{g^2k^2}{4}
(
q^2r_{ik}-v_iv_k
)
-\fr12 gq
(r_{ik}-\fr1{q^2}v_iv_k)
(\wt k_my^m)
+
y^my^na_{nikm}
$$
(which confirms (4.8)),
together with
\ses
$$
\fr B{K^2}y^mb^nR_{nikm}=
y^mb^na_{nikm}
+\fr g{qB}
(S^2b_i-bu_i)u_ja_n{}^j{}_{km}y^mb^n
$$
\ses
$$
=
v_i\wt k_k-(\wt k_jy^j)r_{ik}
+\fr g{qB}
(S^2b_i-bu_i)
[q^2\wt k_k-(\wt k_jy^j)v_k]
$$
and
$$
\fr B{K^2}b^my^nR_{nikm}=
-\fr12   gq(r_{ik}-\fr1{q^2}v_iv_k)
(\wt k_mb^m)
+v_k\wt k_i-(\wt k_jy^j)r_{ik}
$$
supplemented by
$$
\fr B{K^2}b^mb^nR_{nikm}=
b^mb^na_{nikm}
+\fr g{qB}
(S^2b_i-bu_i)u_ja_n{}^j{}_{km}
b^mb^n
$$
\ses
$$
=
-(\wt k_jb^j)r_{ik}
-\fr g{qB}
(S^2b_i-bu_i)
v_k(\wt k_jb^j).
$$


The result of calculations is the following:
$$
g^{nm}R_{nikm}=
\fr{g^2k^2}{4}\Biggl[
(N-2)
\Bigl[
r_{ik}
+\fr g{Bq}(S^2b_i-bu_i)
v_k
\Bigr]
+\fr{ g}{qB}(b+gq)(q^2r_{ik}-v_iv_k)
\Biggr]
$$
\ses
\ses
$$
+\fr{ g}{2q}\Bigl[
-(r_{i}{}^m-\fr1{q^2}v_iv^m)\wt k_mv_k
-v_i(\wt k_k(N-2)-\wt k_mr^m{}_k)\Bigr]
-\fr12  \fr gq
r_{ik}
(\wt k_mv^m)
$$
\ses
$$
-\fr {g^2}{2B}
(S^2b_i-bu_i)\Bigl[
\wt k_k(N-2)-(r_k{}^m-\fr1{q^2}v_kv^m)\wt k_m
\Bigr]
$$
\ses\ses
$$
-
\fr{ g^2}{2B}(b+gq)
(r_{ik}-\fr1{q^2}v_iv_k)
(\wt k_my^m)
$$
\ses\ses
$$
-\fr gq\Biggl[
v_i\wt k_k-(\wt k_jy^j)r_{ik}
+\fr g{qB}
(S^2b_i-bu_i)
[q^2\wt k_k-(\wt k_jy^j)v_k]
\Biggr]
$$
\ses
$$
-\fr gq\Biggl[
-\fr12   gq(r_{ik}-\fr1{q^2}v_iv_k)
(\wt k_mb^m)
+v_k\wt k_i-(\wt k_jy^j)r_{ik}
\Biggr]
$$
\ses
\ses
$$
+\fr{ gb}q\Biggl[
-(\wt k_jb^j)r_{ik}
-\fr g{qB}
(S^2b_i-bu_i)
v_k(\wt k_jb^j)
\Biggr]
$$
\ses\ses\ses\ses
$$
+
a^{nm}a_{nikm}
+\fr g{qB}
\Bigl[(S^2b_i-bu_i)y^ja_{njkm}a^{nm}
+
(b+gq)
y^my^na_{nikm}
\Bigr].
$$


Simplifying yields
$$
g^{nm}R_{nikm}=
\fr{g^2k^2}{4}\Biggl[
(N-2)
\Bigl[
r_{ik}
+\fr g{Bq}(S^2b_i-bu_i)
v_k
\Bigr]
+\fr{ g}{qB}(b+gq)(q^2r_{ik}-v_iv_k)
\Biggr]
$$
\ses
\ses
$$
+\fr{ g}{2q}\Bigl[
-(r_{i}{}^m-\fr1{q^2}v_iv^m)\wt k_mv_k
-v_i(\wt k_kN-\wt k_mr^m{}_k)\Bigr]
+\fr12  \fr gq
r_{ik}
(\wt k_mv^m)
$$
\ses
$$
-\fr {g^2}{2B}
(S^2b_i-bu_i)\Bigl[
\wt k_kN-(r_k{}^m+\fr1{q^2}v_kv^m)\wt k_m
\Bigr]
$$
\ses\ses
$$
-
\fr{ g^2}{2B}(b+gq)
(r_{ik}-\fr1{q^2}v_iv_k)
(\wt k_my^m)
$$
\ses\ses
$$
+\fr gq
(\wt k_jy^j)r_{ik}
+\fr12g^2(r_{ik}-\fr1{q^2}v_iv_k)
(\wt k_mb^m)
-\fr gqv_k\wt k_i
$$
\ses\ses
\be
+
a^{nm}a_{nikm}
+\fr g{qB}
\Bigl[(S^2b_i-bu_i)y^ja_{njkm}a^{nm}
+
(b+gq)
y^my^na_{nikm}
\Bigr].
\ee


Whenever $
 \wt k_n=fb_n,
 $
we obtain
$$
g^{nm}R_{nikm}=
\fr{g^2k^2}{4}\Biggl[
(N-2)
\Bigl[
r_{ik}
+\fr g{Bq}(S^2b_i-bu_i)
v_k
\Bigr]
+\fr{ g}{qB}(b+gq)(q^2r_{ik}-v_iv_k)
\Biggr]
$$
\ses
\ses
$$
-\fr{ gN}{2q}
v_ib_kf
-\fr {g^2}{2B}
(S^2b_i-bu_i)
b_kNf
-
\fr{ g^2}{2B}b(b+gq)
(r_{ik}-\fr1{q^2}v_iv_k)
f
$$
\ses\ses
$$
+\fr{ gb}q
r_{ik}f
+\fr12g^2(r_{ik}-\fr1{q^2}v_iv_k)
f
-\fr gqv_kb_if
$$
\ses\ses
\be
+
a^{nm}a_{nikm}
+\fr g{qB}
\Bigl[(S^2b_i-bu_i)y^ja_{njkm}a^{nm}
+
(b+gq)
y^my^na_{nikm}
\Bigr],
\ee
\ses\ses\ses\ses
\be
R_i{}^h{}_{hk}=
\fr{g^2k^2}{4}
(N-2)r_{ik}
- \fr {gN}{2q}v_ib_kf
+a_i{}^h{}_{hk},
\ee
and
\be
R^{nm}{}_{mn}=
\fr1{K^2}\Biggl[
\fr{g^2k^2}{4}
(N-2)\Bigl((N-1)B+g^2q^2\Bigr)
+\fr12g^2q^2Nf\Biggr]+
a_n{}^h{}_{hm}g^{nm}
\ee
with
\be
a_n{}^h{}_{hm}g^{nm}=
\fr B{K^2}\Bigl(a_n{}^h{}_{hm}a^{nm}
+\fr {gb}q(N-1)f
\Bigr)
+\fr1{K^2}\fr gq(b+gq)y^na_n{}^h{}_{hm}y^m.
\ee


\setcounter{sctn}{5}
\setcounter{equation}{0}

\bc
{\large 5. Warped Product Structure of Associated Riemannian Space}

\ec

\ses\ses

Remarkably, the Landsberg--characteristic condition  (1.15)--(1.16), as well as the Berwald-characteristic
condition  $  \nabla_ib_j=0  $,
does not involve the Finsleroid charge
parameter $g$ and, therefore, imposes restrictions on
only  the  underlying associated Riemannian space
$\cR_N=(M,{\cal S})$.
What is the geometrical meaning of the restrictions?

 The fact that the Finsleroid--axis 1-form $b$ is closed (see (1.15)) and the $b$-lines are geodesics (see
(1.23) and Proposition 1)
can be used to introduce  on the background Riemannian space
$\cR_N=(M,{\cal S})$ the {\it $b$-geodesic coordinates}, to be denoted as $z^i=\{z^0,z^a\}$,
such that with respect to such coordinates we shall have
\be
b^0(z^i)=1, \qquad b^a(z^i)=0, \qquad b_0(z^i)=1, \qquad b_a(z^i)=0,
\ee
\be
a_{00}(z^i)=1, \qquad a_{0a}(z^i)=0, \qquad a_{ab}(z^i)=r_{ab}(z^i), \qquad r_{00}(z^i)=r_{0a}(z^i)=0
\ee
(we have used the unit length (2.5) and the nullification (2.7)),
obtaining  the square of the  Riemannian line element $ds$ of
$\cR_N=(M,{\cal S})$ to be the sum
\be
(ds)^2=(dz^0)^2+r_{ab}(z^i)dz^adz^b
\ee
and the Landsberg--characteristic condition (1.16) to read
\be
\fr12\D{r_{ab}(z^i)}{z^0}=k(z^i)r_{ab}(z^i);
\ee
the indices $a,b,c$ are specified over the range $1,\dots,N-1$.

Solving the last equation by the help of the  representation
\be
r_{ab}(z^i)=\phi^2(z^i)p_{ab}(z^c)
\ee
subject to
\be
\D {\phi}{z^0}=k,
\ee
we can convert the sum (5.3) into the sum
\be
(ds)^2=(dz^0)^2+\phi^2(z^i)p_{ab}(z^c)dz^adz^b
\ee
which says us that $ds$ is  {\it a warped product metric}.
To emphasize the fact that the product property (5.7) has appeared upon
specifying the  $z^0$-coordinate line
to be tangent to the Finsleroid--axis
vector field $b^i(x)$, we shall attribute to the space and metric obtained the quality  of being
of the {\it $b$-warped product type}.
Thus we have the following.

\ses\ses
{\large Theorem 6}. {An  \it $\cF\cF^{PD}_g$--space is of the Landsberg type if and only if
the associated
Riemannian space is a $b$-warped product space and the Finsleroid charge $g$ is a constant.

}
\ses\ses


The warped product type of the Riemannian space was nicely exposed in
the section 13.3 of the book [7] (and we follow the terminology used therein).
In that source, the  interesting special class of warped products that is
characterized by the  condition that the function $\phi(z^i)$ be
 independent of $z^c$ was introduced and considered.
Accordingly, we call the associated Riemannian space {\it special $b$-warped space} if the sum (5.7)
is simplified to be
\be
(ds)^2=(dt)^2+\phi^2(t)p_{ab}(z^c)dz^adz^b,
\ee
where $t=z^0$.
In this case from (5.4)--(5.6) it follows that the function $k$, and hence the field
$\wt k_n=\D k{x^n}+k^2b_n$ introduced in accordance with (1.28), may depend on only $t=z^0$. Therefore,
 when treating with respect to the $b$-geodesic coordinates, we have the proportionality of $\wt k_n$
 to $b_n$. However, the last proportionality is tensorial and must remain valid in terms of any admissible coordinates.
Thus the following theorem may be stated to hold.

\ses\ses
{\large Theorem 7}. {An  \it $\cF\cF^{PD}_g$--space of the Landsberg type is of the
special $b$-warped product type
if and only if the {\rm1}-form $\wt k_n(x)y^n$ is proportional to the {\rm 1}-form $b=b_n(x)y^n$}:
\be
 \wt k_n(x)=f(x)b_n(x).
\ee

\ses\ses

It is the condition (5.9) under which we have obtained the simplified representations
(4.27)--(4.30) in the previous section.



It is also interesting to specialize  the Landsberg--type $\cF\cF^{PD}_g$--space farther
by the help of the following definition:
The Landsbergian $\cF\cF^{PD}_g$--space is called $b${\it-stationary} if
\be
 b_ia_n{}^i{}_{km}=0
 \ee
 holds at any  point $x\in M$.
In this case, given a  $\cF\cF^{PD}_g$--space of the $b$-stationary  Landsberg type
provided the space is not of the Berwald type,
 the identities
\be
A_iR^i{}_k=0, \qquad A_iR^i{}_{km}=0
\ee
hold,
$k$ is not a constant,  the equalities
\be
b_n=\fr{\partial}{\partial x^n}\fr1k, \qquad \wt k_n=0
\ee
are fulfilled,
and the curvature tensor (4.3) is entirely independent of tangent vectors $y$.
With respect to the $b$-stationary coordinates, from (5.1) and (5.12) we obtain $k=1/t$.
If we specify the Landsberg type of the
$\cF\cF^{PD}_g$--space  by stipulating that
\be
(\nabla_m\nabla_n
-\nabla_n\nabla_m)b_k=0,
\ee
then, noting (1.26), we may re-interpret
this commutator   in terms of the Finsler--nature covariant derivatives
to read
\be
b_{k|n|m}-b_{k|m|n}=0,
\ee
so that
the identity (5.10)
together with
 \be
 b_iR_n{}^i{}_{km}=0
 \ee
and
\be
b_iR^i{}_{km}=0, \qquad b_iR^i{}_{k}=0
\ee
(see (1.27))
must appear.
Since  the vanishing
$
y_iR^i{}_{km}=y_iR^i{}_k=0
$
takes place in any Finsler space,
and  $A_i$ are linear combinations of $y_i$ and $b_i$ (see (A.7) in Appendix A),
the  contraction of the tensors
$R^i{}_k$ and $R^i{}_{km}$ by $A_i$ must produce zeros,
whence we arrive at (5.11).
Finally, in view of (1.27), the stipulation (5.13) entails the vanishing (5.12).
The opposite way, namely derivation of the formulas
(5.11)--(5.16) from the stationarity condition (5.10), can easily be
gone, too.
Each of the above formulas
(5.11)--(5.16)  can be taken as the condition that may characterize the $b$-stationary Landsberg case instead of the
departure condition (5.10).


\setcounter{sctn}{6}
\setcounter{equation}{0}

\bc
{\large 6. Indefinite (Relativistic) Case}
\ec

\ses\ses

The positive--definite case (to which the previous four sections were devoted)
can be juxtaposed with  a
 indefinite pseudo--Finsleroid case, assuming the input metric tensor $\{a_{ij}(x)\}$ to be
 pseudo--Riemannian with
the {\it time--space signature:}
\be
{\rm sign}(a_{ij})=(+ - -\dots).
\ee
Namely, attempting  to generalize the pseudo--Riemannian geometry in a pseudo--Finsleroid Finslerian way,
we are to adapt the consideration to the  following decomposition of the tangent bundle $TM$:
\be
TM=\cS_g^+\cup \Si_g^+\cup{\cR_g}\cup\Si_g^-\cup\cS_g^-,
\ee
which sectors relate to the cases that the tangent vectors $y\in TM$ are, respectively,
time--like, upper--cone isotropic, space--like, lower--cone isotropic, or past--like.
The sectors  are defined according to the following list:
\be
\cS_g^+=
\Bigl(y\in \cS_g^+:~y\in T_xM,\,  b(x,y)>-g_-(x)q(x,y)\Bigr),
\ee
\medskip
\be
\Si_g^+=\Bigl(y\in \Si_g^+:~y\in T_xM, \, b(x,y)=-g_-(x)q(x,y)\Bigr),
\ee
\medskip
\be
\cR_g^+=
\Bigl(y\in \cR_g^+:~y\in T_xM,\,  -g_-(x)q(x,y)>b(x,y) >0 \Bigr),
\ee
\medskip
\be
\cR^0=
\Bigl(y\in \cR^0:~y\in T_xM,\,  b(x,y) =0 \Bigr),
\ee
\medskip
\be
\cR_g^-=
\Bigl(y\in \cR_g^-:~y\in T_xM,\,0>b(x,y) >-g_+(x)q(x,y) \Bigr),
\ee
\medskip
\be
\Si_g^-=\Bigl(y\in \Si_g^-:~y\in T_xM, \, b(x,y)=-g_+(x)q(x,y)\Bigr),
\ee
\medskip
\be
\cS_g^-=
\Bigl(y\in \cS_g^-:~y\in T_xM,\,  b(x,y)<-g_+(x)q(x,y)\Bigr),
\ee
\medskip
\be
\cR_g=\cR_g^+\cup\cR_g^-\cup\cR^0.
\ee
We use the convenient notation
\be
G= \fr gh,
\qquad h = \sqrt{1+\fr14g^2}
\ee
(instead of (2.10)),
\be
g_+=-\fr12g+h, \qquad g_-=-\fr12g-h,
\ee
\medskip
\be
G_+=\fr{g_+}h\equiv -\fr12G+1, \qquad G_-=\fr{g_-}h\equiv -\fr12G-1,
\ee
\medskip
\be
g^+=\fr 1{g_+}=-g_-,  \qquad  g^-=\fr 1{g_-}=-g_+,
\ee
\medskip
\be
g^+=\fr12g+h, \qquad g^-=\fr12g-h,
\ee
\medskip
\be
G^+=\fr{g^+}h\equiv \fr12G+1, \qquad G^-=\fr{g^-}h\equiv \fr12G-1.
\ee
The following identities hold
\be
g_++g_-=-g, \qquad g_+-g_-=2h,
\ee
\medskip
\be
g^++g^-=g, \qquad g^+-g^-=2h,
\ee
\medskip
\be
g_+g_-=-1, \qquad
g^+g^-=-1,
\ee
together with
the {\it
$g$--symmetry}
\be
g_+\g -g_-, \qquad g^+\g -g^-, \qquad G_+\g -G_-, \qquad G^+\g -G^-.
\ee

It is implied that $g=g(x)$ is a scalar on the underlying manifold $M$.
All the range
\be
-\infty <g(x)<\infty
\ee
(instead of (2.9)) is now
admissible.
We also  assume that the manifold $M$ admits
a  normalized 1-form $b=b(x,y)$ which is
{\it timelike} in terms of the pseudo--Riemannian metric ${\cal S}$,
 such that the pseudo--Riemannian length of the involved vector $b_i$   be equal to 1.
With respect to  natural local coordinates in the space
$\cR_N$
we have the local representations
\be
 b=b_i(x)y^i,
\ee
\ses
\be
S= \sqrt{|a_{ij}(x)y^iy^j|},
\ee
\ses
\be
q=\sqrt{|r_{ij}(x)y^iy^j|},
\ee
\ses
\be
r_{ij}(x)=b_i(x)b_j(x)-a_{ij}(x),
\ee
\ses
\be
a^{ij}b_ib_j=1,
\ee
\ses
\be S^2=b^2-q^2,
\ee
\ses
\be
b^ir_{ij}=0,
\ee
where
\be
b^i~:=a^{ik}b_k
\ee
(compare with (2.1)--(2.8)).


The {\it  pseudo--Finsleroid characteristic
quadratic form}
\be
B(x,y) :=b^2-gqb-q^2
\equiv (b+g_+q)(b+g_-q)
\ee
is now of the positive discriminant
\be
D_{\{B\}}=4h^2>0
\ee
(compare with (2.11) and (2.12)).

In terms of these concepts, we propose

\ses\ses

{\large Definition}. The scalar function $F(x,y)$ given by the formulas
\be
F(x,y)~:=\sqrt{|B(x,y)|}\,J(x,y)
\equiv
|b+g_-q|^{G_+/2}|b+g_+q|^{-G_-/2}
\ee
and
\be
J(x,y)=
\lf|
\fr{b+g_-q}{b+g_+q}
\rg|^{-G/4},
\ee
is called
the {\it  pseudo--Finsleroid--Finsler  metric function}.

\ses\ses

The positive  (not absolute) homogeneity  holds: $F(x,\la y)=\la F(x,y)$ for any $\la >0$.
The functions
 \be
 L(x,y) =q-\fr g2b
\ee
and
\be
A(x,y)=b-\fr g2q
\ee
are now to be used instead of (2.18) and (2.19), so that (2.20) changes to read
\be
L^2-h^2b^2=B, \qquad A^2-h^2q^2=B.
\ee

Similarly to (2.21), we introduce

\ses\ses

 {\large Definition}.  The arisen  space
\be
\cF\cF^{SR}_g :=\{\cR_{N};\,b(x,y);\,g(x);\,F(x,y)\}
\ee
is called the
 {\it pseudo--Finsleroid--Finsler space}.

\ses\ses

The upperscript ``SR" emphasizes the Specially  Relativistic character of the space under study.

\ses\ses

 {\large Definition}. The space $\cR_N=(M,\,{\cal S})$ entering the above definition is called the {\it associated
pseudo--Riemannian space}.

\ses\ses

 {\large Definition}. The scalar $g(x)$ is called
the {\it  pseudo--Finsleroid charge}.
The 1-form $b$ is called the  {\it  pseudo--Finsleroid--axis}  1-{\it form}.

\ses\ses

It can be verified that the Finslerian metric tensor constructed from the function $F$ given by (6.32)
does inherit from the tensor $\{a_{ij}(x)\}$
the time--space signature (6.1):
\be
{\rm sign}(g_{ij})=(+ - -\dots).
\ee

The structure (2.27) for the
 Cartan tensor remains valid in the pseudo--Finsleroid case, now with
 \be
A_hA^h=-\fr{N^2}{4}g^2.
\ee
Elucidating the structure of
the respective indicatrix curvature tensor
 (2.29) of the  $\cF\cF^{SR}_g$--space again results in the special type (2.30),  with
$S^*=\fr14g^2$,
so that
$$
 \cR_{\text{pseudo-Finsleroid  Indicatrix} }=-(1+\fr14g^2) \le -1
$$
(compare with (2.32)) and
$$
\cR_{\text{pseudo-Finsleroid Indicatrix} }\stackrel{g\to 0}{\Longrightarrow}
\cR_{\text{pseudo-Euclidean Sphere}}=-1.
$$
{\it The pseudo--Finsleroid indicatrix
 is a space of constant  negative curvature}.

By analogy with
(2.33) and (2.34)
we obtain the {\it $\cF\cF_g^{SR}$--scalar product}
 \be
<y_1,y_2>_{\{x\}}~:=F(x,y_1)F(x,y_2) \cosh\Bigl( \al_{\{x\}}(y_1,y_2)\Bigr),\qquad y_1,y_2\in\cS_g^+,
 \ee
 and the $\cF\cF_g^{SR}$--{\it angle}
\ses
\be
 \al_{\{x\}}(y_1,y_2)~: = \fr1h\arccosh \fr{ A(x,y_1)A(x,y_2)-h^2<y_1,y_2>_{\{x\}}^{\{r\}} }
 {\sqrt{|B(x,y_1)|}\,\sqrt{|B(x,y_2)|} }, \qquad y_1,y_2\in\cS_g^+,
 \ee
where again $<y_1,y_2>_{\{x\}}^{\{r\}}=r_{ij}(x)y_1^iy_2^j$.
The $\cF\cF_g^{SR}$--{\it distance} $\De s$ between ends of the vectors in the tangent space
is now given by the formula
\be
 (\Delta s)^2 = (F(x,y_1))^2 + (F(x,y_2))^2 - 2 F(x,y_1)F(x,y_2) \cosh\Bigl( \al_{\{x\}}(y_1,y_2)\Bigr).
 \ee
Also,
\be
 <y,y>_{\{x\}}=(F(x,y))^2.
  \ee

The pseudo--Finsleroid analogue of the angle (2.39) between  the Finsleroid axis and a vector
now reads
 \be
 \al_{\{x\}}(y)~: = \fr1h\arccosh \fr{ A(x,y)}
 {\sqrt{|B(x,y)|}}.
 \ee

\ses

All the  motivation presented  in  previous Sections 4 and 5
 when elucidating the Berwald and Landsberg cases,
can be repeated  word--for--word
in the present indefinite type,
with  slight changes
in  formulas.

\ses\ses

\setcounter{sctn}{7}
\setcounter{equation}{0}

\bc
{\large 7. Conclusions}
\ec

\ses\ses

Thus, we have  geometrized the tangent bundle
by the help of the
Finsleroid, {\it resp.} pseudo--Finsleroid, in the positive--definite approach, {\it resp.} indefinite approach
of the relativistic signature.
Under the particular condition $g=const$,
the Finsleroid metric function $K$, {\it resp.}  the pseudo--Finsleroid metric function $F$,
is defined by a 1-form $b$ and a Riemannian, {\it resp.}  pseudo--Riemannian, metric. In general, they involve also
a scalar, the ``geometric charge`` $g(x)$.

Both the metric functions, $K$ defined by (2.14) and $F$ defined by (6.32),
are positively, and {\it not } absolutely, homogeneous of the degree 1:
$K(x,\la y) =\la  K(x,y),\, \la>0$, together with
$F(x,\la y) =\la  F(x,y),\, \la>0$.
Also, they are reversible in the extended sense
$$
K_{\bigl|g\to -g, ~y\to -y}=K, \qquad F_{\bigl|g\to -g, ~y\to -y}=F.
$$
The traditional symmetry assumption $K(x,-y)= K(x,y)$, as well as $F(x,-y)= F(x,y)$, remains valid for the vectors
$y\in T_xM$ obeying $b(x,y)=0$, and is violating otherwise unless $g=0$.

Each  metric function, $K$ and $F$,
is  rotund around the Finsleroid axis and, therefore, reflects
the idea of the
 spherical symmetry in the spaces of directions orthogonal to the axis.
Their indicatrices are spaces of constant curvature, namely positive in case of $K$ and negative in case of $F$.

The Finslerian and pseudo--Finsleroid spaces arisen are naturally
endowed with the inner scalar product and the concept of angle.

We have deduced the affirmative and detailed
answers to all the two questions posed in the beginning of the present paper.
Moreover, we have supplemented the answers by explicit and rather simple representations
for the key geometric objects associated thereto.
A particular  peculiarity is that the respective  Landsberg--case tensor
$A_{i|j}$
 proves to be proportional to the lengthened angular metric tensor
${\cal H}_{ij}$, according to (1.19).

Surprisingly, the tensor ${\cal H}_i{}^j$
proves to be entirely independent of the Finsleroid charge $g$
(see (A.16) in Appendix A below), that is
of whether the Finsleroid extension of the underlying Riemannian
geometry is performed or not. In this respect, the tensor
is {\it universal} and is free tractable in, and coming from,  the entire context of the underlying
 Riemannian space.


  Having an $\cF\cF^{PD}_g$--space with constant Finsleroid charge,
calculations of the $hv$--curvature tensor $P_k{}^i{}_{mn}$
yield the   result explicitly in terms of the lengthened angular metric tensor
${\cal H}_{mn}$ (see (1.18) and the end of Appendix A below)
and, in the Landsberg--case thereto,
the structure of the associated $hh$-curvature tensor $R_n{}^i{}_{km}$
can clearly be explained
by the help of the explicit representations (4.1)--(4.30) obtained, in particular
the Ricci--deflection vector and tensor are
remarkably proportional to the tensor
${\cal H}_n{}^m$ in accordance with (4.21) and (4.22).

The comprehension of the
Landsberg case  $\cF\cF^{PD}_g$--space in terms of the category of a warped product structure
(to which Section 5 was devoted)
seems to be a constructive and handy idea.


The $C^{\infty}$-class smoothness, though   appears frequently in publications devoted to Finsler geometry
and might outwardly sound as being a desirable attractive tune,
is ordinarily ill-incorporating in designs of particular Finslerian metric functions.
Our examples $K$ and $F$ used are not  total exclusions from this rule,
namely the functions are analytic along any direction except for that which points in, or opposite to, the direction
of the Finsleroid axis, that is, when  $q\to 0$.
This actually does not destroy the potential  applied abilities of the Finsleroid metric functions,
as long as we exercise a sufficient care when encountering with the Finsleroid--axis--directions.
The metric functions $K$ and $F$
and all the components of the associated covariant vector $y_i$,
as well as the determinant (see (A.6) in Appendix A below)
of the  Finslerian metric tensor,
are all everywhere smooth, the divergences (at $q=0$)  start arising with  components
of the  Finslerian metric tensor.
In particular,
 $g_{ij}$ and  $g^{ij}$ (see (A.2) and (A.3) in Appendix A below), as well as
 $A_i$ and $A^i$ (see (A.7) and (A.8) in Appendix A below),
  are $(\sim 1/q)$-singular,
and simultaneously  the contraction $A_hA^h$
(see (2.28)) is entirely independent of $y$.
The $hv$-curvature tensor in the Landsberg case is also singular when $q\to 0$
(see (1.18)), however the very condition (1.16) that underlines the Landsberg case
is free of any singularities. On the other hand, the components of
the tensor $q^2{\cal H}_i{}^j$ are quadratic in the variables $y$ and, therefore, are  analytical in any direction
(examine the right--hand side of the representation (A.16) placed in Appendix A below).

By following the book [7], we say that the Finsler space is {\it $y$-global} if the space is smooth and strongly
convex on all of the {\it slit tangent bundle }
$TM\smallsetminus 0$. In many cases it is important to
 supplement the notion by explicit indicating the degree of smoothness. The following assertions are valid.
The $\cF\cF^{PD}_g $--space  is  $y$-global of the class  $C^1$, and not of the class $C^2$,
  on all of the $TM\smallsetminus 0$.
The $\cF\cF^{PD}_g $--space  is  $y$-global of the class  $C^{\infty}$
on all of the {\it $b$-slit  tangent bundle }
\be
{\cal T}_bM~:= TM \smallsetminus0\smallsetminus b\smallsetminus -b
\ee
(obtained
by deleting out in $TM\smallsetminus 0$
all the directions which point along, or oppose, the directions  given rise to by the Finsleroid--axis 1-form $b$;
the value $q=0$ just corresponds to such directions in view of the initial formulas (2.3) and (2.7)),
provided
the  $C^{\infty}$--smoothness is assumed for  the input Riemannian space
$\cR_N=(M,S)$, the 1-form $b$,
and the scalar
$g(x)$.


Above our consideration was everywhere referred to local vicinities of the points $x\in M$ of the underlying manifold.
Let us ponder over the possibility to construct global
$\cF\cF^{PD}_g $--structures.
 In the Berwald case, that is when $k=0$, the tensor
$r_{ab}(z^i)$ looses the dependence on $z^0$ (as this is evident from  (5.4)), so that the underlying
Riemannian space $\cR_N=(M,{\cal S})$ must be locally the Cartesian product of
 an $(N-1)$-dimensional Riemannian space $\cR_{N-1}$ and  an interval of
 a  one--dimensional Euclidean space
 $\set R$, with
  the $r_{ab}(z^c)$ playing the role of the Riemannian metric tensor of the former space
and the Finsleroid--axis 1--form $b$ being tangent to the latter space.
  In particular, on the product Riemannian spaces
\be
\cR_N=\cR_{N-1}\times \set R \quad {\text{and}} \quad  \cR_N=\cR_{N-1}\times \set S^1,
\ee
  the Berwald--case
$\cF\cF^{PD}_g $--space structure can be introduced straightforwardly by imposing on such  spaces
the Finsleroid metric function $K$ in accordance with its  very definition (1.14) and assuming the
Finsleroid charge $g$ to be a constant. If we modify the
construction in but
multiplying the metric tensor
 $r_{ab}(z^c)$ of the involved space
   $\cR_{N-1}$ by a  factor
$\phi(z^i)$,
 we just lift the structure from the Berwaldian level to the
 Landsbergian  type, with the function $k$ obtained from (5.6).
The Landsberg--type  $\cF\cF^{PD}_g $--structures
 thus constructed on the grounds   of  the product Riemannian  spaces (7.2) are $y$-global of the class $C^1$
(and not of the class $C^2$) on the slit tangent bundle $TM\smallsetminus 0$,
  and are  $y$-global of the class $C^{\infty}$ on the $b$-slit tangent bundle (7.1)
  (assuming the input Riemannian space
     $\cR_{N-1}$ and  the Finsleroid--axis 1-form $b$  to be of the class $C^{\infty}$).


{\it In the two--dimensional case}, $N=2$, we have identically ${\cal H}_{ij}=0$.
 Therefore, the implication $(3.1)~\to~(3.2)$
noted in Section 3 reduces to $Landsberg~type ~ \to~Berwald~type$ and the representation (1.10)  reduces to merely
$\dot A_i=0$, which says us the following:
Any  two--dimensional   $\cF\cF^{PD}_g $--space is of the Berwald type,
 whenever the Finsleroid charge $g(x)$ is a constant
via the condition (1.6), the Riemannian case being occurred if only $g=0$.
The Landsberg scalar $J$ is zero if and only if the Finsleroid charge $g(x)$ is a constant.
 The two--dimensional analog to the Cartan tensor structure (2.27) is
the representation $A_{ijk}=I(x)A_iA_jA_k$ with the main scalar $I(x)=|g(x)|$ which is the function of position $x$ only.
The components
of $A_i$ are singular at only $q=0$ (the formula (A.7)  in Appendix A is applicable to the two--dimensional case), whence
the quantity $I(x)$ is  meaningful on all the indicatrix $I_x\subset T_xM$ except for the
top and down points.
At first sight, a suspicion may arise that we contradict to the known and remarkable theorem
(see pp. 278 and 279 in the book [7]) which claims that the independence of the main scalar of tangent vectors
entails strictly the Riemannian case, that is, $I=I(x)$ entails $I=0$. However,
an attentive consideration of the conditions under which the theorem has been settled out
just reveals the fact that the conditions imply  the $y$-global
smoothness of at least class $C^3$ to be valid  for the Finsler space.
It is the lack of the $y$-global property of the class $C^3$  that is  the reason
why the two--dimensional  $\cF\cF^{PD}_g $--space is spared the fate of
degenerating to become  Riemannian and occurs being of the Berwald type proper.
All the formulae  obtained in Section 4 for the $hh$-curvature tensor are applicable at N=2.
In particular, the formula (4.6) shows that $\it Ric$ at $N=2$ is not zero.


{\it Comparison with the Randers space} brings to evidence many interesting interrelations. Let the dimension $N$
be more than 2.
It should be noted first of all that the nature of the Berwald case can easily be clarified in
the  $\cF\cF^{PD}_g $--space as well as in the Randers space, resulting in quite a similar condition
for the characteristic  1-form $b=b_k(x)y^k$
to be  parallel with respect to the associated Riemannian metric, $\nabla_i b_j=0$.
In both  the $\cF\cF^{PD}_g $--space and the Randers space, the weakly Berwald  condition
${ A}_{j|i}=0$  entails the complete Berwald condition
${ A}_{jmn|i}=0$.
There is, however,  a significant distinction in the structures of the Randers space and the
 $\cF\cF^{PD}_g $--space in respect of the Landsberg--type properties.  Namely,
 the identical vanishing
$
A^j{\dot A}_j=0
$
may be considered to be the weakest necessary condition to underline the full Landsberg condition
$\dot A_{ikj}=0$.
It proves that  the identical vanishing
of the contraction
$
A^j{\dot A}_j
$ is {\it always} valid  in the   $\cF\cF^{PD}_g $--space with a constant Finsleroid
charge and {\it never }  in the non--Riemannian
Randers space. The latter fact says us  that the Randers space can {\it never}
be of the properly Landsberg type, that is, for the Randers space
the notions of being Lansberg and being Berwald are equivalent.
Can the condition $\nabla_i b_j=k(x)(a_{ij}-b_ib_j)$ generate a non--Berwaldian Landsberg space?
 The true answer, namely  ``yes'' in the  $\cF\cF^{PD}_g $--space
{\it versus} ``no'' in the Randers space,
 can be motivated in simple words. Namely, in both spaces the dotted vector $\dot A_i$ under the condition
 $\nabla_i b_j=k(x)(a_{ij}-b_ib_j)$,
as well as the very field $A_i$, is spanned by a linear
combination of two vectors, $y_i$ and $b_i$, so that such two scalars $c_1,c_2$ must exist that
$\dot A_i=c_1y_i+c_2A_i$. Since $y^i\dot A_i$ is  zero in any
Finsler space, we must put $c_1=0$. The watershed is that in  $\cF\cF^{PD}_g $--space, -- and not in the Randers space, --
we also have  zero for the contraction $A^i\dot A_i$ which entails that  $c_2$ is also zero.
So in  $\cF\cF^{PD}_g $--space we can safely  pose the condition
$\nabla_i b_j=k(x)(a_{ij}-b_ib_j)$ to get a non--Berwaldian Landsberg type,
while such a door is closed for ever in  the Randers space.
Obviously, the distinction is rooted in difference of the algebraic structures of the Cartan tensors associated with
the  $\cF\cF^{PD}_g $--space and with the Randers space.
It is also interesting to note for comparison
that the condition of the type $\nabla_i b_j=\fr12\si (a_{ij}-b_ib_j)$  with
$\si=const$ appears when characterizing the Randers spaces of constant flag curvature
(see [10,11]).


\setcounter{equation}{0}

\bc
{\large  Appendix A:  Representations for Distinguished $\cF\cF^{PD}_g $--Objects }
\ec

\ses\ses

With the notation $u_k=y^na_{nk}$,
the direct calculations of the covariant tangent vector and Finslerian metric tensor
on the basis of the Finsleroid metric function $K$ given by (2.14)  yield
\be
y_i=(a_{ij}y^j+gqb_i)\fr{K^2}B
\ee
and
\be
g_{ij}=
\biggl[a_{ij}
+\fr g{B}\Bigl((gq^2-\fr{bS^2}q)b_ib_j-\fr bqu_iu_j+
\fr{ S^2}q(b_iu_j+b_ju_i)\Bigr)\biggr]\fr{K^2}B.
\ee
The reciprocal components $(g^{ij})=(g_{ij})^{-1}$ read
\be
g^{ij}=
\biggl[a^{ij}+\fr gq(bb^ib^j-b^iy^j-b^jy^i)+\fr g{Bq}(b+gq)y^iy^j
\biggr]\fr B{K^2}.
\ee
We also obtain
\be
y_ib^i=(b+gq)\fr{K^2}B, \qquad
g_{ij}b^j=\bigl(b_i+gq\fr{y_i}{K^2}\bigr)\fr{K^2}B,
\ee
and
\be
h_{ij}b^j=\bigl(b_i-b\fr{y_i}{K^2}\bigr)\fr{K^2}B.
\ee
The determinant of the tensor (A.2) is everywhere positive:
\be
\det(g_{ij})=\bigl(\fr{K^2}B\bigr)^N\det(a_{ij})>0.
\ee


By the help of the formulas (A.2)--(A.4) we find
\be
A_i=\fr {NK}2g\fr1{q}(b_i-\fr b{K^2}y_i)
\ee
and
\be
A^i=\fr N2g\fr 1{qK}
\Bigl[Bb^i-(b+gq)y^i\Bigr]
\ee
together with
\be
A_ib^i=\fr N2gq\fr KB, \qquad    A^ib_i=\fr N2gq\fr 1K.
\ee
These formulas are convenient to verify the algebraic structure (2.27) and the contraction (2.28).
Also,
\be
A_{ij}~:= K\partial A_i/\partial y^j+l_iA_j=
-\fr N2\fr{gb}{q}{\cal H}_{ij}+\fr{2}NA_iA_j
\ee
with
\be
{\cal H}_{ij}=h_{ij}-\fr{A_iA_j}{A_nA^n}.
\ee


The last tensor fulfills obviously the identities
\be
{\cal H}_{ij}y^j=0,\qquad
{\cal H}_{ij}A^j=0,
\ee
which in turn entails
\be
{\cal H}_{ij}b^j=b_j{\cal H}_{i}{}^{j}=0
\ee
because $A^i$ are linear combinations of $y^i$ and $b^i$
(see (A.8)).
We also have
\be
g^{ij}{\cal H}_{ij}=N-2,
\ee
\ses
\be
 g^{mn}{\cal H}_{im}{\cal
H}_{jn}={\cal H}_{ij},
\ee
\ses
\be
{\cal H}_i{}^j=\de_i{}^j-b_ib^j-\fr1{q^2}(a_{in}y^n-bb_i)(y^j-bb^j)
\equiv r_i{}^j-\fr1{q^2}v_iv^j,
\ee
and
\be
K\Bigl(\D{{\cal H}_{ij}}{y^k}
-\D{{\cal H}_{kj}}{y^i}
\Bigr)
=l_k{\cal H}_{ij}-l_i{\cal H}_{kj}
-\fr1{A_nA^n}
\fr{Ngb}{2q}
(A_k{\cal H}_{ij}-A_i{\cal H}_{kj}
)
\ee
together with
$$
\D{(q^2{\cal H}_i^j)}{y^k}
-\D{(q^2{\cal H}_k^j)}{y^i}
=3\Bigl[
(\de_i{}^j-b_ib^j)(a_{kn}y^n-bb_k)
-
(\de_k{}^j-b_kb^j)(a_{in}y^n-bb_i)
\Bigr]
$$
\ses
\be
=3\Bigl[
{\cal H}_i^j(a_{kn}y^n-bb_k)
-
{\cal H}_k^j(a_{in}y^n-bb_i)
\Bigr].
\ee

The  structure (2.27) of the $\cF\cF^{PD}_g $--space Cartan tensor is such that
\be
A_k\3Aikj=
\fr1N(A_iA_j+h_{ij}A_kA^k)=\fr1N(2A_iA_j+{\cal H}_{ij}A_kA^k),
\ee
so that the tensor
\be
\tau_{ij}~:
=
A_{ij}-A_k\3Aikj
\ee
is equal to
\be
\tau_{ij}=-\fr N4\fr{g(2b+gq)}{q}{\cal H}_{ij}.
\ee
By comparing (A.21) with  (A.12) and (A.13), the identities
$
\tau_{ij}y^j=0,~\tau_{ij}A^j=0,
$
and
\be
 \tau_{ij}b^j=b_j\tau_i{}^j=0
\ee
just follow.


The tensor
\be
\tau_{ijmn}~:=K\partial A_{jmn}/\partial y^i-A_{ij}{}^hA_{hmn}-A_{im}{}^hA_{hjn}
-A_{in}{}^hA_{hjm}
+l_jA_{imn}+l_mA_{ijn}+l_nA_{ijm}
\ee
can be expressed as follows:
\be
\tau_{ijmn}=-\fr{g(2b+gq)}{4q}({\cal H}_{ij}{\cal H}_{mn}+{\cal H}_{im}{\cal H}_{jn}+{\cal H}_{in}{\cal H}_{jm}),
\ee
showing the total symmetry in all four indices and the properties
$$
\tau_{ij}= g^{mn}\tau_{ijmn}
$$
and
\be
\qquad y^i\tau_{ijmn}=0,\qquad A^i\tau_{ijmn}=0,\qquad b^i\tau_{ijmn}=0.
\ee

We use
 the Riemannian covariant derivative
\be
\nabla_ib_j~:=\partial_ib_j-b_ka^k{}_{ij},
\ee
where
\be
a^k{}_{ij}~:=\fr12a^{kn}(\prtl_ja_{ni}+\prtl_ia_{nj}-\prtl_na_{ji})
\ee
are the
Christoffel symbols given rise to by the associated Riemannian metric ${\cal S}$.

For   the $hh$-curvature tensor $R^i{}_k$ we use the formula
\be
K^2R^i{}_k~:=
2\D{\bar G^i}{x^k}-\D{\bar G^i}{y^j}\D{\bar G^j}{y^k}
-y^j\Dd{\bar G^i}{x^j}{y^k}
+2\bar G^j\Dd{\bar G^i}{y^k}{y^j}
\ee
(which is tantamount to the definition (3.8.7) on p. 66 of the book [7];
$\bar G^i=
\fr12\ga^i{}_{nm}y^ny^m$,
with the Finslerian Christoffel symbols
$
\ga^i{}_{nm}$).
The concomitant tensors
\be
R^i{}_{km}~:=\fr1{3K}\Biggl(\D{(K^2R^i{}_{k})}{y^m}-\D{(K^2R^i{}_{m})}{y^k}
\Biggr)
\ee
and
\be
R_n{}^i{}_{km}~:=\D{(KR^i{}_{km})}{y^n}
\ee
arise.


The cyclic identity
\be
R_j{}^i{}_{kl|t}+R_j{}^i{}_{lt|k}+R_j{}^i{}_{tk|l}=
P_j{}^i{}_{ku}R^u{}_{lt}+P_j{}^i{}_{lu}R^u{}_{tk}+P_j{}^i{}_{tu}R^u{}_{kl}
\ee
(the formula (3.5.3) on p. 58 of the book [7]) is valid in any Finsler space.
If we contract (A.31) by $g^{jl}\de^k{}_i$, we get:
\be
g^{jl}\Bigl(R_j{}^i{}_{il|t}+R_j{}^i{}_{lt|i}+R_j{}^i{}_{ti|l}\Bigr)=
P^l{}^i{}_{iu}R^u{}_{lt}+P^l{}^i{}_{lu}R^u{}_{ti}+P^l{}^i{}_{tu}R^u{}_{il}.
\ee
Under the Landsberg condition, the tensor $P_{ijkl}$ is symmetric in all of its four indices and, therefore,
the previous identity reduces to merely
\be
g^{jl}\Bigl(R_j{}^i{}_{il|t}+R_j{}^i{}_{lt|i}+R_j{}^i{}_{ti|l}\Bigr)=0,
\ee
which can also be written as
\be
\rho^i{}_{j|i}\equiv 0
\ee
with
\be
\rho_{ij}~:=\fr12(R_i{}^m{}_{mj}+R^m{}_{ijm})-\fr12g_{ij}R^{mn}{}_{nm}.
\ee
This conservation law (A.34)--(A.35) is valid in any Finsler space of the Landsberg type and, therefore,
in any Landsberg case of the $\cF\cF^{PD}_g $--space studied.

\ses

To verify (1.31)--(1.32) the reader can use
the  relationship
\be
{\dot A}_{jkl}= -\fr14y_i
\fr{\partial^3\ga^i{}_{nm}y^ny^m}
{\partial y^j\partial y^k\partial y^l}
\ee
(the formula (3.8.5)
on p. 65 of the book [7])
 which is valid in any Finsler space.
It is easy to check that the insertion of
the particular coefficients (1.24)
in (A.36)
results in  the Landsbergian
${\dot A}_{jkl}=0$.


Finally, we are able to compute  components of the $hv$-curvature tensor $P_k{}^i{}_{mn}$.
To this end we use the formula
\be
P_k{}^i{}_{mn}=K\fr{\partial\Bigl({\dot A}^i{}_{km}-\fr12G^i{}_{km}\Bigr)}{\partial y^n}
\ee
(see (3.8.4) on p. 65 of [7]) which is valid in any Finsler space. In the $\cF\cF^{PD}_g $--space
under the only condition that the Finsleroid charge
is a constant,  direct tedious
calculations yield the   result explicitly in terms of the lengthened angular metric tensor
${\cal H}_{mn}$, namely
\ses\ses\\
$$
\fr 2KP_k{}^i{}_{mn}=
\fr gq
\Biggl(f^i{}_k
\mu_{mn}
+f^i{}_m
\mu{}_{kn}
-(\nabla_kb_m+\nabla_mb_k)
\mu^i{}_{n}
\Biggr)
$$
\ses\ses
$$
+\fr g{q^3}\Biggr(
t_k\mu^i{}_{mn}-q^2(\nabla_kb_n)
\mu^i{}_{m}
+t_m\mu^i{}_{kn}-q^2(\nabla_mb_n)
\mu^i{}_{k}
-t^i\mu_{kmn}+
q^2(\nabla^ib_n)
\mu_{km}
\Biggr)
$$
\ses\ses
$$
-\fr {g^2}{2q^2}\Biggl[
\Biggl(b^j\nabla_jb^i
-\fr1{q^2}v^iy^hb^j\nabla_jb_h
\Biggr)
\mu_{kmn}
-\Biggl(
v_nb^j\nabla_jb^i
-v^ib^j\nabla_jb_n
\Biggr)
\mu_{km}
\Biggr]
$$
\ses\ses
$$
-
\fr {g^2}{2q^2}\Biggr[
\Biggl(
b^j\nabla_jb_k
-\fr1{q^2}v_ky^hb^j\nabla_jb_h
\Biggr)
\mu^i{}_{mn}
-
\Biggl(
v_nb^j\nabla_jb_k
-v_kb^j\nabla_jb_n
\Biggr)
\mu^i{}_m
\Biggr]
$$
\ses\ses
$$
-
\fr {g^2}{2q^2}\Biggr[
\Biggl(
b^j\nabla_jb_m
-\fr1{q^2}v_my^hb^j\nabla_jb_h
\Biggr)
\mu^i{}_{kn}
-
\Biggl(
v_nb^j\nabla_jb_m
-v_mb^j\nabla_jb_n
\Biggr)
\mu^i{}_k
\Biggr]
$$
\ses\ses
\be
-\fr{g^2}{2q^2}y^hb^j\nabla_jb_h
\Biggr(
\mu^i{}_n
\mu_{km}
+
\mu^i{}_m
\mu_{kn}
+
\mu^i{}_k
\mu_{mn}
\Biggr),
\ee
\ses\\
where
\ses\ses\ses\\
\be
\mu^i{}_n=r^i{}_n-\fr1{q^2}v^iv_n\equiv {\cal H}^i{}_{n}, \qquad \mu_{in}=r_{in}-\fr1{q^2}v_iv_n
\equiv \mu^j{}_na_{ji}\equiv
{\cal H}_{in}\fr B{K^2},
\ee
\ses
\be
\mu^i{}_{kn}
=
r^i{}_kv_n+r^i{}_nv_k+r_{nk}v^i-\fr3{q^2}v^iv_nv_k\equiv
{\cal H}^i{}_{k}v_n+{\cal H}^i{}_{n}v_k+{\cal H}_{nk}\fr B{K^2}v^i,
\ee
\ses
$$
\mu_{kmn}=
r_{km}v_n+r_{kn}v_m+r_{mn}v_k-\fr3{q^2}v_kv_mv_n
\equiv
\mu^j{}_{mn}a_{jk}
$$
\ses
\be\equiv
\Bigl({\cal H}_{km}v_n+{\cal H}_{kn}v_m+{\cal H}_{mn}v_k\Bigr)\fr B{K^2},
\ee
\ses\ses\\
and
\be
t_k=y^h\nabla_kb_h,   \qquad t^i=a^{ik}t_k.
\ee
It is easy to observe that
\be
P_k{}^i{}_{mn}y^n=0, \qquad P_k{}^i{}_{mn}=P_m{}^i{}_{kn},
\ee
and
\be
y^kP_k{}^i{}_{mn}=-\dot A^i{}_{mn},
\ee
where the right--hand side agrees exactly with  the right--hand side of
(1.31).
It is also easy to verify that upon the Landsbergian condition (1.16)
the above tensor $P_k{}^i{}_{mn}$ given by (A.38) reduces to the tensor
(1.18).

\ses
\ses

\end{document}